\documentclass[twoside,11pt]{article}
\usepackage{journal2e_v1}

\usepackage[numbers]{natbib}
\usepackage{epsf}
\usepackage{epsfig}
\usepackage{amsmath}
\usepackage{subfigure}
\usepackage{amssymb}
\usepackage{algorithm}
\usepackage{multirow}
\usepackage{multicol}

\usepackage{graphicx} 

\usepackage{algorithm}
\usepackage{algorithmic}

\usepackage{epsf}
\usepackage{epsfig}
\usepackage{amsmath}
\usepackage{subfigure}
\usepackage{amssymb}
\usepackage{multirow}
\usepackage{multicol}
\usepackage{pgfplots}
\usepackage{algorithm}
\usepackage{algorithmic}

\usepackage{enumitem}

\usepackage{color} 

\usepackage{hyperref}       
\usepackage{url}            

\def\0{{\bf 0}}
\def\1{{\bf 1}}

\def\OM{{\mathcal O}}

\def\SM{{\mathcal S}}

\def\RB{{\mathbb R}}

\def\EB{{\mathbb E}}
\def\PB{{\mathbb P}}

\newcommand{\ti}[1]{\tilde{#1}}

\def\argmin{\mathop{\rm argmin}}

\def\diag{\mathrm{diag}}

\def\RB{{\mathbb R}}
\def\tV{{\tilde{V}}}
\def\tH{{\tilde{H}}}
\def\hV{{\hat{V}}}
\def\hL{{\hat{L}}}

\def\tI{{\tilde{I}}}

\newcommand{\norm}[1]{\left\|#1\right\|}
\title{Approximate Newton Methods}
\author{ \name Haishan Ye \\
	\addr {hsye\_cs@outlook.com}\\
	Shenzhen Research Institute of Big Data\\ 
	The Chinese University of Hong Kong, Shenzhen
	\AND 
	\name Luo Luo   \\
\addr {rickyluoluo@gmail.com } \\
Department of Computer Science and Engineering \\
Shanghai Jiao Tong University \\
800 Dong Chuan Road, Shanghai, China 200240	
\AND
\name Zhihua Zhang   \\
\addr zhzhang@math.pku.edu.cn \\
School of Mathematical Sciences \\
Peking  University \\
Beijing, China 100871
}
\begin{document}

\maketitle

\begin{abstract}%
	Many machine learning models involve solving optimization problems. Thus, it is important to deal with a large-scale optimization problem in big data applications. 
	Recently, subsampled Newton methods have emerged to attract much attention 
	due to their efficiency at
	each iteration, rectified a weakness in the ordinary Newton method of suffering a high cost in each
	iteration while commanding a high convergence rate. Other efficient stochastic second order methods are also proposed. However, the convergence properties of these methods are still not well understood. There are also several important gaps between the current convergence theory and the performance in real applications. In this paper, we aim to fill these gaps. We propose a unifying framework to analyze both local and global convergence properties of second order methods. Based on this framework, we present our theoretical results which match the performance in real applications well. 
\end{abstract}

\section{Introduction}

Mathematical optimization is an important pillar of machine learning. We consider the following optimization problem:
\begin{equation}
\min_{x\in\RB^{d}} F(x) = \frac{1}{n}\sum_{i=1}^{n}f_i(x), \label{eq:prob_desc}
\end{equation}
where the $f_i(x)$ are smooth functions. 
Many machine learning models can be expressed as \eqref{eq:prob_desc} where each $f_i$ is the loss with respect to (w.r.t.) the $i$-th training sample. There are many examples such as logistic regressions, smoothed support vector machines, neural networks, and graphical models. 

Many optimization algorithms to solve the problem in~\eqref{eq:prob_desc} are based on the following iteration:
\[
x^{(t+1)} = x^{(t)} - s_t Q_t \text{g}(x^{(t)}), \; t=0, 1, 2, \ldots,
\]
where $s_t>0$ is the step length.
If $Q_t$ is the identity matrix and $\text{g}(x^{(t)}) = \nabla F(x^{(t)})$, the resulting procedure is called \emph{Gradient Descent} (GD) which achieves sublinear convergence for a general smooth convex objective function and linear convergence for a smooth-strongly convex objective function. When $n$ is large, the full gradient method is inefficient due to its iteration cost scaling linearly in $n$. Consequently, stochastic gradient descent (SGD) has been a typical alternative \cite{robbins1951stochastic,li2014efficient,cotter2011better}. In order to achieve cheaper cost in each iteration, such a method constructs an approximate gradient on a small mini-batch of data. However, the convergence rate can be significantly slower than that of the full gradient methods \cite{nemirovski2009robust}. Thus, a great deal of efforts have been made to devise modification to achieve the convergence rate of the full gradient while keeping low iteration cost \cite{johnson2013accelerating, roux2012stochastic,schmidt2013minimizing,Zhang}.

If $Q_t$ is a $d\times d$ positive definite matrix of containing the curvature information, this formulation leads us to \emph{second-order} methods. It is well known that second order methods enjoy superior convergence rate in both theory and practice in contrast to \emph{first-order} methods which only make use of the gradient information. The standard Newton method, where $Q_t = [\nabla^2 F(x^{(t)})]^{-1}$,  $\text{g}(x^{(t)}) = \nabla F(x^{(t)})$ and $s_t = 1$,  achieves a  quadratic  convergence rate for smooth-strongly convex objective functions. However, the \emph{Newton method} takes $\OM(nd^2+d^3)$ cost per iteration, so it becomes extremely expensive when $n$ or $d$ is very large. As a result, one tries to construct an approximation of the Hessian in which way the update is computationally feasible,  while keeping sufficient second order information. One class of such methods are quasi-Newton methods, which are generalizations of the secant methods to find the root of the first derivative for multidimensional problems. The celebrated Broyden-Fletcher-Goldfarb-Shanno (BFGS) and its limited memory version (L-BFGS) are the most popular and widely used \cite{nocedal2006numerical}. They take $\OM(nd+d^2)$ cost per iteration.

Recently,  when $n\gg d$, a class of called \emph{subsampled Newton} methods have been proposed, which define an approximate Hessian matrix with a small subset of samples.  The most naive approach is to sample a subset of functions $f_i$ randomly \cite{roosta2019sub,byrd2011use,xu2016sub} to construct a subsampled Hessian. \citet{erdogdu2015convergence} proposed  a  regularized subsampled Newton method called NewSamp. When the Hessian can be written as $\nabla^2F(x) = [B(x)]^T B(x)$ where $B(x)$ is an available $n\times d$ matrix, \citet{pilanci2015newton} used sketching techniques to approximate the Hessian and proposed \emph{sketch Newton} method. Similarly, \citet{xu2016sub} proposed to sample rows of $B(x)$ with non-uniform probability distribution. \citet{agarwal2016second} brought up an algorithm called LiSSA to approximate the inverse of Hessian directly.

Although the convergence performance of stochastic second order methods has been analyzed, the convergence properties are still not well understood. There are several important gaps lying between the convergence theory and the performance of these algorithms in real applications.

First, it is about  the necessity of Lipschitz continuity of the Hessian. In previous work, to achieve a linear-quadratic convergence rate, stochastic second order methods all assume that $\nabla^2F(x)$ is Lipschitz continuous. However, in real application without this assumption, they might also converge to an optimal point. For example, \citet{erdogdu2015convergence} used NewSamp to  successfully
train the smoothed-SVM  in which the $\ell_2$-hinge loss is used, so the corresponding Hessian is not Lipschitz continuous.

Second, it involves the sketched size of sketch Newton methods. To obtain a linear convergence, the sketched size is $\OM(d\kappa^2)$ in \cite{pilanci2015newton} and then be improved to $\OM(d\kappa)$ in \cite{xu2016sub} using Gaussian sketching matrices, where $\kappa$ is the condition number of the Hessian matrix in question. However, the sketch Newton empirically performs well even when the Hessian matrix is ill-conditioned. Sketched size being several tens of times, or even several times of $d$ can achieve a linear convergence rate in unconstrained optimization. But the theoretical result of \citet{pilanci2015newton,xu2016sub}  implies that sketched size may be beyond $n$ in ill-condition cases.

Third,  it talks about the sample size in regularized subsampled Newton methods.  In both \cite{erdogdu2015convergence} and \cite{roosta2016sub}, their theoretical analysis shows that the sample size of regularized subsampled Newton methods should be set as the same as the conventional subsampled Newton method. In practice, however, adding a large regularizer can obviously reduce the sample size while keeping convergence.  Thus, this does not agree with  the extant theoretical analysis \cite{erdogdu2015convergence,roosta2016sub}. 


In this paper, we aim to fill these gaps between the current theory and empirical performance. More specifically, we first  cast these second order methods into an algorithmic  framework that we call \emph{approximate Newton}. Accordingly,  we propose a general result for analysis of both
local and global convergence properties of second order methods. 
Based on this framework, we then give detailed theoretical analysis which matches the empirical performance. We summarize our contribution as follows:	
\begin{itemize}
	\item We propose a unifying framework (Theorem~\ref{thm:univ_frm} and Theorem~\ref{thm:glb}) to analyze local and global convergence properties of second order methods including stochastic  and deterministic versions. 
	The convergence performance of second order methods can be analyzed easily and systematically in this framework. 
	\item We prove that the Lipschitz continuity condition of Hessian is not necessary for achieving linear and superlinear convergence in variants of subsampled Newton. But it is needed to obtain quadratic convergence. This explains the phenomenon that NewSamp \cite{erdogdu2015convergence} can be used to train the smoothed SVM in which the Lipschitz continuity condition of Hessian is not satisfied. It also reveals the reason why previous stochastic second order methods, such as subsampled Newton, sketch Newton, LiSSA, etc., all achieve a linear-quadratic convergence rate. 
	\item We prove that the sketched size is \emph{independent} of the condition number of Hessian matrix which explains that sketched Newton performs well even when Hessian matrix is ill-conditioned. 
	\item Based on our analysis framework, we provide a much tighter bound of sample size of subsampled Newton methods. 
	To the best knowledge of authors, it is the tightest bound of subsampled Newton methods.
	\item We provide a theoretical guarantee that adding a regularizer is an effective way to reduce sample size in subsampled Newton methods while keeping converging. Our theoretical analysis also shows that adding a regularizer will lead to poor convergence behavior as the sample size decreases.
\end{itemize}

\subsection{Organization}	

The remainder of the paper is organized as follows.	In Section~\ref{sec:prelim} we present notation and preliminaries. 
In Section~\ref{sec:frame} we present a unifying framework for local and global convergence analysis of second order methods. In Section~\ref{sec:ske_newton} we analyze the  convergence properties of sketch Newton methods and prove that sketched size is independent of the condition number of Hessian matrix. In Section~\ref{sec:sub_newton} we give the convergence behaviors of several variants of subsampled Newton method. Especially, we reveal the relationship among sample size, regularizer and convergence rate. In Section~\ref{sec:exp}, we validate our theoretical results experimentally. Finally, we conclude our work in Section~\ref{sec:conclusion}.
Theorems are proved in appendices in the order of their appearing.

\section{Notation and Preliminaries}
\label{sec:prelim}	

Section~\ref{subsec:notation} defines the notation used in this paper. Section~\ref{subsec:ske_mat} introduces matrices sketching techniques and their properties. Section~\ref{subsec:ass} describes some important assumptions about objective functions.

\subsection{Notation}\label{subsec:notation}

Given a matrix $A=[a_{ij}] \in \RB^{m \times n}$ of rank $ \ell $ and a positive integer $k\leq \ell$,  its condensed SVD is given as
$A=U\Sigma V^{T}=U_{k} \Sigma_{k} V_{k}^{T}+U_{\setminus k} \Sigma_{{\setminus} k} V_{{\setminus}k}^{T}$,
where $U_{k}$ and $U_{{\setminus}k}$ contain the left singular vectors of $A$,  $V_{k}$ and $V_{{\setminus}k}$ contain the right
singular vectors of $A$, and $\Sigma=\diag(\sigma_1, \ldots, \sigma_{\ell})$ with $\sigma_1\geq \sigma_2 \geq \cdots \geq \sigma_{\ell}>0$ are
the nonzero singular values of $A$. We will use $\sigma_{\max}(A)$ to denote the largest singular value and $\sigma_{\min}(A)$ to denote the smallest non-zero singular value. Thus, the condition number of $A$ is defined by $\kappa(A) \triangleq \frac{\sigma_{\max}(A)}{\sigma_{\min}(A)}$. If $A$ is positive semidefinite, then $U = V$ and the square root of $A$ can be defined as $A^{1/2} = U\Sigma^{1/2}U^T$. It also holds that $\lambda_i(A) = \sigma_i(A)$, where $\lambda_i(A)$ is the $i$-th largest eigenvalue of $A$, $\lambda_{\max}(A) = \sigma_{\max}(A)$, and $\lambda_{\min}(A) = \sigma_{\min}(A)$.

Additionally, $\|A\|_{F} \triangleq (\sum_{i,j}a_{ij}^{2})^{1/2}=(\sum_{i}\sigma_{i}^{2})^{1/2}$
is the Frobenius norm of $A$ and
$\|A\|\triangleq \sigma_{1}$ is the spectral norm. Given a positive definite matrix $M$, $\|x\|_M \triangleq \|M^{1/2}x\|$ is called the $M$-norm of $x$.
Give square matrices $A$ and $B$ with the same size, 	
we	denote $A \preceq B$ if $B-A$ is positive semidefinite. 

\subsection{Randomized sketching matrices} \label{subsec:ske_mat}

We first give an $\epsilon_0$-subspace embedding property which will be used to sketch Hessian matrices. Then we list some useful types of randomized sketching matrices including  Gaussian projection \citep{halko2011finding,johnson1984extensions}, 
leverage score sampling \cite{drineas2006sampling}, count sketch \citep{clarkson2013low, nelson2013osnap,meng2013low}.

\begin{definition} \label{lem:sub-embed}
$S\in\RB^{\ell\times n}$ is said to be an $\epsilon_0$-subspace embedding matrix w.r.t.\ a fixed matrix $A\in\RB^{n\times d}$ where $d<n$, if  $\|SA x\|^{2}=(1\pm\epsilon_0)\|Ax\|^{2}$ \emph{(i.e., $(1-\epsilon_0) \|Ax\|^2 \leq \|SA x\|^{2} \leq (1+\epsilon_0) \|Ax\|^{2}$)} for all $x \in \RB^{d}$.
\end{definition}

From the definition of the $\epsilon_0$-subspace embedding matrix, we can derive the following property directly.
\begin{lemma}
$S\in\RB^{\ell\times n}$ is an $\epsilon_0$-subspace embedding matrix w.r.t.\ the matrix $A\in\RB^{n\times d}$ if and only if
	\[
	(1-\epsilon_0)A^TA\preceq A^TS^TSA\preceq(1+\epsilon_0)A^TA.
	\]
\end{lemma}
\paragraph{Gaussian sketching matrix.}	
The most classical sketching matrix is the Gaussian sketching matrix $S\in\RB^{\ell\times n}$, whose extries are i.i.d.\ from  the normal of mean 0 and variance $1/\ell$. 
Owing to the  well-known concentration properties \cite{woodruff2014sketching},  Gaussian random matrices are very attractive. Besides, $\ell = \OM(d/\epsilon_0^2)$ is enough to guarantee the $\epsilon_0$-subspace embedding property for any fixed matrix $A\in\RB^{n\times d}$. Moreover, $\ell = \OM(d/\epsilon_0^2)$ is the tightest bound among known types of sketching matrices. However, the Gaussian random matrix is usually dense, so it is costly to compute $SA$. 

\paragraph{Leverage score sketching matrix.}
A leverage score sketching matrix $S = D \Omega \in\RB^{\ell\times n}$ w.r.t.\ $A\in\RB^{n\times d}$ is defined by sampling probabilities $p_i$, a sampling matrix $\Omega\in\RB^{n\times \ell}$ and a diagonal  rescaling matrix $D \in\RB^{\ell \times \ell}$. Specifically, we construct $S$ as follows.  For every $j = 1,\dots,\ell$, independently and with replacement, pick an index $i$ from the set $\{1,2\dots,n\}$ with probability  $p_{i}$, and set $\Omega_{ji} = 1$ and $\Omega_{j k}=0$ for $k \neq i$ as well as $D_{jj}=1/\sqrt{p_{i}\ell}$. 
The sampling probabilities $p_i$ are the leverage scores of $A$  defined as follows. Let $V\in\RB^{n\times d}$  be the column orthonormal basis of $A$, and let $v_{i,*}$ denote the $i$-th row of $V$. Then  $q_{i} \triangleq \|v_{i,*}\|^{2}/d$ for $i=1, \ldots, n$ are the leverage scores of $A$. To achieve an $\epsilon_0$-subspace embedding property w.r.t.\ $A$, $\ell = \OM(d\log d/\epsilon_0^2)$ is sufficient. 

\paragraph{Sparse embedding matrix.}	
A sparse embedding matrix $S\in\RB^{\ell\times n}$ is such a matrix in each column of which there is only one nonzero entry uniformly sampled from $\{1,-1\}$ \cite{clarkson2013low}. Hence, it is very efficient to compute $SA$, especially when $A$ is  sparse. To achieve an $\epsilon_0$-subspace embedding property w.r.t.\ $A\in\RB^{n\times d}$, $\ell = \OM(d^2/\epsilon_0^2)$ is sufficient \cite{meng2013low,woodruff2014sketching}.

Other sketching matrices such as Subsampled Randomized Hadamard Transformation \cite{drineas2012fast,halko2011finding} as well as their properties can be found in the survey \cite{woodruff2014sketching}.


\subsection{Assumptions and Notions} \label{subsec:ass}

In this paper, we focus on the problem described in Eqn.~\eqref{eq:prob_desc}. Moreover, we will make the following two assumptions. 
\paragraph{Assumption 1}
The objective function  $F$ is $\mu$-strongly convex, that is, 
\[
F(y) \geq F(x) + [\nabla F(x)]^T(y-x) + \frac{\mu}{2}\|y-x\|^2, \mbox{ for }\; \mu>0. 
\]
\paragraph{Assumption 2}
$\nabla F(x)$ is $L$-Lipschitz continuous, that is, 
\[
\|\nabla F(x) - \nabla F(y)\| \leq L\|y-x\|, \mbox{ for }\; L>0.
\] 
Assumptions 1 and 2 imply that for any $x \in \RB^d$, we have 
\[
\mu I\preceq \nabla^2F(x) \preceq L I, 
\]
where $I$ is the identity matrix of appropriate size. With a little confusion, we define 
\[
\kappa \triangleq \frac{L}{\mu}.
\]
Note that  $\kappa$ is an upper bound of the condition number of the Hessian matrix $\nabla^2F(x)$ for any $x$. 
Furthermore, if $\nabla^2 F(x)$ is Lipschitz continuous, then we have
\[
\|\nabla^2F(x) - \nabla^2F(y)\| \leq \hat{L}\|x-y\|,
\]
where $\hat{L}>0$ is the Lipschitz constant of $\nabla^2F(x)$.

Throughout this paper, we use notions of linear convergence rate, superlinear convergence rate and quadratic convergence rate. In our paper, the convergence rates we will use are defined w.r.t. $\|\cdot\|_{M}$, where $M = \nabla^2F(x^*)$ and $x^*$ is the optimal solution to Problem~\eqref{eq:prob_desc}. A sequence of vectors $\{x^{(t)}\}$ is said to converge linearly to a limit point $x^*$, if for some $0< \rho<1$,
\[
\limsup_{t\to \infty} \frac{\|x^{(t+1)} - x^*\|_{M}}{\|x^{(t)} - x^* \|_{M}} = \rho.
\]
Similarly, superlinear convergence and quadratic convergence are respectively defined as
\begin{align*}
\limsup_{t \to \infty} \frac{\|x^{(t+1)} - x^*\|_{M}}{\|x^{(t)} - x^* \|_{M}} = 0, \qquad
\limsup_{t \to \infty} \frac{\|x^{(t+1)} - x^*\|_{M}}{\|x^{(t)} - x^* \|_{M}^2}= \rho.
\end{align*}
We call it the linear-quadratic convergence rate if the following condition holds: 
\[
\|x^{(t+1)} - x^*\|_{M} \leq \rho_1\|x^{(t)} - x^*\|_{M}+ \rho_2 \|x^{(t)} - x^*\|_{M}^2,
\]
where $0<\rho_1<1$.

\section{Main Results} \label{sec:frame}

The existing variants of stochastic second order methods share some important attributes. First, these methods such as NewSamp \cite{erdogdu2015convergence}, LiSSA \cite{agarwal2016second}, subsampled Newton with conjugate gradient \cite{byrd2011use}, and subsampled Newton with non-uniformly sampling \cite{xu2016sub}, all have the same convergence properties; that is, they have a linear-quadratic convergence rate. 

Second, they also enjoy the same algorithm procedure summarized as follows. In each iteration, they first construct an approximate Hessian matrix $H^{(t)}$ such that
\begin{align}
(1-\epsilon_0)H^{(t)} \preceq\nabla^2F(x^{(t)}) \preceq(1+\epsilon_0)H^{(t)}, \label{eq:prec_cond}
\end{align}  
where $0\leq\epsilon_0<1$. 
Then they solve the following optimization problem
\begin{align}
\min_p \frac{1}{2}p^TH^{(t)}p -p^T\nabla F(x^{(t)}) \label{eq:lin_eq}
\end{align}
approximately or exactly to obtain the direction vector $p^{(t)}$. Finally, their update equation is given as $x^{(t+1)} = x^{(t)} - p^{(t)}$.
With this procedure, we regard these stochastic second order methods as \emph{approximate Newton} methods. The detailed algorithmic description is listed in Algorithm~\ref{alg:app_newton}.

\begin{algorithm}[tb]
\caption{Approximate Newton.}
\label{alg:app_newton}
\begin{small}
	\begin{algorithmic}[1]
		\STATE {\bf Input:} $x^{(0)}$, $0<\delta<1$, $0<\epsilon_0<1$;
		\FOR {$t=0,1,\dots$ until termination}
		\STATE Construct an approximate Hessian $H^{(t)}$ satisfying Condition~\eqref{eq:prec_cond};
		\STATE Calculate $p^{(t)} \approx \argmin_{p} \frac{1}{2} p^T H^{(t)} p - p^T\nabla F(x^{(t)})$;
		\STATE Update $x^{(t+1)}= x^{(t)}-p^{(t)}$;
		\ENDFOR
	\end{algorithmic}
\end{small}
\end{algorithm}

\subsection{Local Convergence Analysis}
In the following theorem, we propose a unifying framework which describes the convergence properties of the second order optimization procedure depicted above.

\begin{theorem}\label{thm:univ_frm}
Let Assumptions 1 and 2 hold. Suppose that $\nabla^2 F(x)$ exists and is continuous in a neighborhood of a minimizer $x^*$. $H^{(t)}$ is a positive definite matrix that satisfies Eqn.~\eqref{eq:prec_cond} with $0\leq\epsilon_0<1$. Let $p^{(t)}$ be an approximate solution of Problem~\eqref{eq:lin_eq} such that
	\begin{equation}
	\|\nabla F(x^{(t)}) - H^{(t)}p^{(t)}\| \leq \frac{\epsilon_1}{\kappa^{3/2}} \|\nabla F(x^{(t)})\|,\label{eq:inexact_p}
	\end{equation} 
	where $0<\epsilon_1<1$. Then Algorithm~\ref{alg:app_newton} has the following convergence properties.
	
	(a) There exists a sufficient small value $\gamma$ and $\nu = o(1)$ such that when $\|x^{(t)} - x^*\|_M\leq \gamma$, we have that 
	\begin{equation}
	\label{eq:lin_conv}
	\small
	\norm{x^{(t+1)} - x^*}_M \leq \left(\epsilon_0 + \epsilon_1 + 2\nu\mu^{-1}+2\left(2\nu^{1/2}\mu^{-1/2} + \nu\mu^{-1}\right)(\nu\mu^{-1}+1)\right)\norm{x^{(t)} - x^*}_M.
	\end{equation}
	Moreover, $\nu$ will go to $0$ as $x^{(t)}$ goes to $x^*$.
	
	(b) Furthermore, if $\nabla^2F(x)$ is $\hat{L}$-Lipschitz continuous, and $x^{(t)}$ satisfies
	\begin{equation}
	\|x^{(t)} - x^{*}\|_M \leq \mu^{3/2}\hL^{-1}, \label{eq:qua_cond}
	\end{equation}
	 then it holds that
	 \begin{equation}
	 \label{eq:lin_qua_conv}
	 \norm{x^{(t+1)} - x^*}_M  
	 \leq 
	 \left(\epsilon_0+\epsilon_1\right) \norm{x^{(t)} - x^*}_M  + 7 \mu^{-3/4} \hL^{1/2} \norm{x^{(t)} - x^*}_M^{3/2}.
	 \end{equation}
\end{theorem}
\begin{remark}
	In Eqn.~\eqref{eq:lin_qua_conv}, the high order term is linear to $\norm{x^{(t)} - x^*}_M^{3/2}$ instead of $\norm{x^{(t)} - x^*}_M^{2} $ in previous work \citep{erdogdu2015convergence,agarwal2016second}.
	However, this difference can be neglected.
	If $\{x^{(t)}\}$ converges with rate $\norm{x^{(t+1)} - x^*}_M \leq \OM\left( \norm{x^{(t)} - x^*}_M^{3/2}\right)$, then it takes $\OM\left(\log_{3/2}\log(1/\epsilon)\right)$ iterations to achieve an $\epsilon$-suboptimality.
	In contrast, 
	if $\{x^{(t)}\}$ converges with rate $\norm{x^{(t+1)} - x^*}_M \leq \OM\left( \norm{x^{(t)} - x^*}_M^{2}\right)$, then it takes  $\OM\left(\log_{2}\log(1/\epsilon)\right)$ iterations.
	Since it holds that $$\log_{3/2}\log{1/\epsilon} = \log_{3/2}2 \cdot\log_{2}\log(1/\epsilon), \;\mbox{and}, \; \log_{3/2}2 < 2,$$
	we will also call the sequence $\{x^{(t)}\}$ satisfying $\norm{x^{(t+1)} - x^*}_M \leq \OM\left( \norm{x^{(t)} - x^*}_M^{3/2}\right)$ converges quadratically . 
	Similarly, we will refer Eqn.~\eqref{eq:lin_qua_conv} as the linear-quadratic convergence.
\end{remark}
From Theorem~\ref{thm:univ_frm}, we can find some important insights. First, Theorem~\ref{thm:univ_frm} provides sufficient conditions to get different convergence rates including linear, and super-liner rates. 
If $\left(\epsilon_0 + \epsilon_1\right)$ is a constant, then sequence $\{x^{(t)}\}$ converges linearly because $\nu = o(1)$ and it will go to $0$ as $t$ goes to infinity.  
Furthermore, if we set $\epsilon_0 = \epsilon_0(t)$ and $\epsilon_1 = \epsilon_1(t)$ such that $\epsilon_0(t)$ and $\epsilon_1(t)$ decrease to $0$ as $t$ increases, then sequence $\{x^{(t)}\}$ will converge super-linearly.

Second, Theorem~\ref{thm:univ_frm} makes it clear  that the Lipschitz continuity of the Hessian is \emph{not necessary} for linear and super-linear convergence  of stochastic second order methods including Subsampled Newton method, Sketch Newton, NewSamp, etc. 
This reveals the reason why NewSamp can be used to train the smoothed SVM where the Lipschitz continuity of the  Hessian matrix  is not satisfied. 
The Lipschitz continuity condition is only needed to get a quadratic convergence or linear-quadratic convergence. 
This explains the phenomena that LiSSA\cite{agarwal2016second}, NewSamp~\cite{erdogdu2015convergence}, subsampled Newton with non-uniformly sampling \cite{xu2016sub}, Sketched Newton \cite{pilanci2015newton} have linear-quadratic convergence rate because they all assume that the Hessian is Lipschitz continuous. 
In fact, it  is well known that the Lipschitz continuity condition of $\nabla^2F(x)$ is not necessary to achieve a linear or superlinear convergence rate for inexact Newton methods.

Third, the unifying framework of Theorem~\ref{thm:univ_frm} contains not only stochastic second order methods, but also the deterministic versions. For example, letting $H^{(t)} = \nabla^2F(x^{(t)})$ and using conjugate gradient to get $p^{(t)}$, we obtain the famous ``Newton-CG'' method. In fact, different choice of $H^{(t)}$ and different way to calculate $p^{(t)}$ lead us to different second order methods. 

\begin{algorithm}[tb]
	\caption{Approximate Newton with backtracking line search.}
	\label{alg:damp_app_newton}
	\begin{small}
		\begin{algorithmic}[1]
			\STATE {\bf Input:} $x^{(0)}$, $0<\alpha<0.5$, $0<\beta<1$;
			\FOR {$t=0,1,\dots$ until termination}
			\STATE Construct an approximate Hessian $H^{(t)}$ satisfying Condition~\eqref{eq:prec_cond};
			\STATE Calculate $p^{(t)} \approx \argmin_{p} \frac{1}{2} p^T H^{(t)} p - p^T\nabla F(x^{(t)})$;
			\STATE Line search: 
				\WHILE {$F(x^{(t)}+s p^{(t)})>F(x^{(t)}) + \alpha s [\nabla F(x^{(t)})]^T p^{(t)}$ }
				\STATE $s = \beta s$
				\ENDWHILE
			\STATE Update $x^{(t+1)}= x^{(t)}-s p^{(t)}$;
			\ENDFOR
		\end{algorithmic}
	\end{small}
\end{algorithm}	

\subsection{Global Convergence Analysis}	

In the previous analysis, the theory is local and approximate Newton can achieve a fast convergence rate once the iterations enter a suitable basin of the origin. 
In this section, we are going to obtain global convergence results for \emph{self-concordant} functions. 
The self-concordant assumption is widely studied in the global convergence analysis of Newton methods \cite{pilanci2015newton,boyd2004convex}. 

Note that a closed, convex function $F$: $\RB^d \to \RB$ is called self-concordant if:
\begin{equation*}
\frac{d}{d\alpha}\nabla^2F(x + \alpha v)|_{\alpha = 0} \preceq 2\|v\|_{x} \nabla^2 F(x)
\end{equation*}
for all $x$ in the domain of $F(x)$ and  and $v \in\RB^d $, where $ \|v\|_x= (v^T \nabla^2F(x)v)^{1/2}$ is the local norm.

To achieve a global convergence, approximate Newton method should combine with the line search.
At the damped phase where $[\nabla F(x^{(t)})]^T p^{(t)}$ is large, line search is applied to guarantee the convergence of  approximate Newton method.
Once $[\nabla F(x^{(t)})]^T p^{(t)}$ is sufficient small, then step size $s = 1$ can keep approximate Newton converging with a linear rate.
The detailed algorithmic description of approximate Newton with backtracking line search is listed in Algorithm~\ref{alg:damp_app_newton}.

In the following theorem, we provide the iteration complexity of  Algorithm~\ref{alg:damp_app_newton} to achieve an $\epsilon$-suboptimality.
\begin{theorem}
	\label{thm:glb}
	Assuming the objective function $F(x)$ is self-concordant, $H^{(t)}$ is a positive definite matrix satisfying Eqn.~\eqref{eq:prec_cond} with $0\leq\epsilon_0<1$. 
	Let $p^{(t)}$ be a descent direction  satisfying Eqn.~\eqref{eq:inexact_p}. 
	The total complexity of approximate Newton method with backtracking line search (Algorithm~\ref{alg:damp_app_newton}) to achieve an $\epsilon$-suboptimality   is at most 
	\begin{equation}
	T = \frac{F(x^{(0)}) - F(x^*)}{\eta} + \frac{2}{1-\epsilon_0- 2\epsilon_1\kappa^{-1} }\log\left(\frac{1-\epsilon_0 - 2\epsilon_1\kappa^{-1}}{12\epsilon}\right), \label{eq:iter_comp}
	\end{equation} 
	where $\eta$ is defined as
	\begin{equation*}
	\eta = \alpha\beta \frac{(1-\epsilon_0)\rho^2(1-\epsilon_0 - 2\epsilon_1\kappa^{-1})^2}{144 + 12 \rho \sqrt{(1-\epsilon_0)} (1-\epsilon_0 - 2\epsilon_1\kappa^{-1})}, \quad\mbox{with}\quad 	
	\rho =\frac{\left(1-\epsilon_1\kappa^{-1}\cdot\left(\frac{1+\epsilon_0}{1-\epsilon_0}\right)^{1/2}\right)^{1/2} }{(1+\epsilon_0)^{1/2} \left(1+\epsilon_1\kappa^{-1}\cdot\left(\frac{1+\epsilon_0}{1-\epsilon_0}\right)^{1/2}\right)}.
	\end{equation*}
\end{theorem}
\begin{remark}
	In the above theorem, the iteration complexity of approximate Newton with line search still depends the condition number of the objective function even it is self-concordant. 
	This dependence on the condition number is caused by the approximation to $H^{-1}\nabla F(x)$. 
	If $\epsilon_1 = 0$ in Eqn.~\eqref{eq:inexact_p}, then we can obtain that $\eta = \alpha\beta\frac{(1-\epsilon_0)^3}{144(1+\epsilon_0) + 12 (1+\epsilon_0)^{1/2}(1-\epsilon_0)^{3/2}}$ which is independent of the condition number.
	Thus, the total complexity is independent of the independent of the condition number.
\end{remark}

\begin{algorithm}[tb]
	\caption{Sketch Newton.}
	\label{alg:sketch_newton}
	\begin{small}
		\begin{algorithmic}[1]
			\STATE {\bf Input:} $x^{(0)}$, $0<\delta<1$, $0<\epsilon_0<1$;
			\FOR {$t=0,1,\dots$ until termination}
			\STATE Construct an $\epsilon_0$-subspace embedding matrix $S$ for $B(x^{(t)})$ and where $\nabla^2 F(x)$ is of the form $ \nabla^2 F(x) = (B(x^{(t)}))^TB(x^{(t)})$, and calculate $H^{(t)} = [B(x^{(t)})]^{T}S^{T}SB(x^{(t)})$;
			\STATE Calculate $p^{(t)} \approx \argmin_{p} \frac{1}{2} p^T H^{(t)} p - p^T\nabla F(x^{(t)})$;
			\STATE Update $x^{(t+1)}= x^{(t)}-p^{(t)}$;
			\ENDFOR
		\end{algorithmic}
	\end{small}
\end{algorithm}

\section{Sketch Newton Method}\label{sec:ske_newton}

In this section, we use Theorem~\ref{thm:univ_frm} to analyze  the convergence properties of Sketch Newton which utilizes the sketching technique to approximate the Hessian. We mainly focus on the case that the Hessian matrix is of the form
\begin{align}
\nabla^2 F(x) = B(x)^TB(x) \label{eq:H_ass}
\end{align}
where $B(x)$ is an explicitly available $n\times d$ matrix. Our result can be easily extended to the case that 
\begin{align*}
\nabla^2 F(x) = B(x)^TB(x) + Q(x),
\end{align*}	  
where $Q(x)$ is a positive semi-definite matrix related to the Hessian of regularizer. 

The Sketch Newton method constructs the approximate Hessian matrix as follows:
\begin{equation}
\label{eq:ske_H}
H^{(t)} =[S^{(t)}B(x)]^T S^{(t)}B(x)
\end{equation}  
where $S^{(t)} \in \RB^{\ell \times n}$ is a randomized sketching matrix. 
Approximate Newton method with such Hessian approximation is referred as  sketch Newton method. 
The detailed algorithmic description is listed in Algorithm~\ref{alg:sketch_newton}.
\begin{theorem}\label{thm:Sketch_newton}
	Let $F(x)$ satisfy the conditions described in Theorem~\ref{thm:univ_frm}. Assume the Hessian matrix is given as Eqn.~\eqref{eq:H_ass}. Let $0<\delta<1$, $0<\epsilon_0<1/2$ and $0\leq \epsilon_1 <1$ be given. $S\in\RB^{\ell \times n}$ is an $\epsilon_0$-subspace embedding matrix w.r.t.\ $B(x)$ with probability at least $1-\delta$. Then sketch Newton (Algorithm~\ref{alg:sketch_newton})  has the following convergence properties:	
	\begin{enumerate}[label = (\alph*)]
		\item There exists a sufficient small value $\gamma$ and $\nu = o(1)$ such that when $\|x^{(t)} - x^*\|_M\leq \gamma$,  each iteration satisfies  Eqn.~\eqref{eq:lin_conv} with probability at least $1-\delta$.
		\item If $\nabla^{2}F(x^{(t)})$ is also Lipschitz continuous and $\{x^{(t)}\}$ satisfies Eqn.~\eqref{eq:qua_cond}, then each iteration satisfies Eqn.~\eqref{eq:lin_qua_conv} with probability at least $1-\delta$.
		\item  If $F(x)$ is furthermore self-concordant, the iteration complexity of the sketch Newton with backtracking line search (Algorithm~\ref{alg:damp_app_newton} with $H^{(t)}$ constructed as Eqn.~\eqref{eq:ske_H}) is upper bounded by Eqn.~\eqref{eq:iter_comp}.
	\end{enumerate}
\end{theorem}

\begin{table}[]
	\centering
	\caption{Comparison with previous work}
	\label{tb:comp}
	\bgroup
	\def\arraystretch{1.5}	
	\begin{tabular}{ccc}
		\hline
		Reference & Sketched Size & Condition number free? \\ \hline
		\citet{pilanci2015newton} &     $\OM\left(\frac{d\kappa^2\log d}{\epsilon_0^2}\right)$               &    No               \\ \hline
		\citet{xu2016sub}   &  $\OM\left(\frac{d\kappa\log d}{\epsilon_0^2}\right)$                  &     No               \\ \hline
		{\bf Our result} (Theorem~\ref{thm:Sketch_newton})&   $\OM\left(\frac{d\log d}{\epsilon_0^2}\right)$                 &     Yes              \\ \hline
	\end{tabular}
	\egroup
\end{table}	

Theorem~\ref{thm:Sketch_newton} directly provides a bound of the sketched size. Using the leverage score sketching matrix as an example, the sketched size $\ell = \OM(d\log d/\epsilon_0^2)$ is sufficient. We compare our theoretical bound of the sketched size with the ones of \citet{pilanci2015newton} and \citet{xu2016sub} in Table~\ref{tb:comp}. As we can see, our sketched size is much smaller than the other two, especially when the Hessian matrix is ill-conditioned.
Theorem~\ref{thm:Sketch_newton} shows that the sketched size $\ell$ is \emph{independent} on the condition number of the Hessian matrix $\nabla^2F(x)$ just as shown in Table~\ref{tb:comp}. This explains the phenomena that when the Hessian matrix is ill-conditioned, Sketch Newton performs well even when the sketched size is only several times of $d$.

Furthermore, the iteration complexity of the sketch Newton with backtracking line search shares the similar result to the one of \citet{pilanci2015newton}.
Especially when $\epsilon_1 = 0$, Eqn.~\eqref{eq:iter_comp} reduces to 
\begin{equation*}
T = \frac{F(x^{(0)}) - F(x^*)}{\eta} + 4\log\left(\frac{1}{24\epsilon}\right),\;\; \mbox{with}\quad \eta = \frac{\alpha\beta(1-\epsilon_0)^3}{12\left((1+\epsilon_0) +  (1+\epsilon_0)^{1/2}(1-\epsilon_0)^{3/2}\right)}.
\end{equation*} 
We can observe that $T$ is independent of the condition number of the objective function. 
A similar result can be found in Theorem 2 of \citet{pilanci2015newton}.

Theorem~\ref{thm:Sketch_newton} also contains the possibility of achieving an asymptotically super-linear rate by using an iteration-dependent sketching accuracy $\epsilon_0 = \epsilon_0(t)$. In particular, we present the following corollary.

\begin{corollary} \label{cor:suplin}
$F(x)$ satisfies the the properties described in Theorem~\ref{thm:univ_frm}. Consider the approximate Hessian $H^{(t)}$ constructed as Eqn.~\eqref{eq:ske_H} with the iteration-dependent sketching accuracy is given as $\epsilon_0(t) = \frac{1}{\log(1+t)}$ and $p^{(t)} = [H^{(t)}]^{-1}\nabla F(x)$. If the initial point $x^{(0)}$ is close enough to the optimal point $x^*$, then sequence $\{x^{(t)}\}$ of the sketch Newton (Algorithm~\ref{alg:app_newton} with $H^{(t)}$ constructed as Eqn.~\eqref{eq:ske_H}) converges superlinearly.
\end{corollary}

\section{The Subsampled Newton method and Variants}\label{sec:sub_newton}	

In this section, we apply Theorem~\ref{thm:univ_frm} to analyze  subsampled Newton methods. 
Instead of the Hessian can be presented as Eqn.~\eqref{eq:H_ass}, for subsample Newton methods, we assume that the Hessian be the sum of different Hessian's:
\begin{equation}
\nabla^2 F(x) = \frac{1}{n} \sum_{i=1}^{n}\nabla^2 f_i(x), \quad\mbox{with}\quad \nabla^2 f_i(x) \in\RB^{d\times d}.
\end{equation}
We make the assumption that each $f_i(x)$ and $F(x)$ have the following properties:
\begin{align}
\max_{1 \leq i\leq n}\|\nabla^2f_i(x)\| \leq K < \infty,\label{eq:k} \\
\lambda_{\min}(\nabla^2F(x))\geq \mu>0. \label{eq:sigma}
\end{align}
Accordingly, we can define a new kind of condition number $\hat{\kappa} = \frac{K}{\mu}$.

\begin{algorithm}[tb]
	\caption{Subsampled Newton.}
	\label{alg:H_subsamp}
	\begin{small}
		\begin{algorithmic}[1]
			\STATE {\bf Input:} $x^{(0)}$, $0<\delta<1$, $0<\epsilon_0<1$;
			\STATE Set the sample size $|\mathcal{S}| \geq \frac{3K/\mu\log(2d/\delta)}{\epsilon_0^2}$.
			\FOR {$t=0,1,\dots$ until termination}
			\STATE Select a sample set $\SM$ of size $|\SM|$, and $H^{(t)} = \frac{1}{|\SM|}\sum_{j\in\mathcal{S}}\nabla^2 f_j(x^{(t)})$;
			\STATE Calculate $p^{(t)} \approx \argmin_{p} \frac{1}{2} p^T H^{(t)} p - p^T\nabla F(x^{(t)})$;
			\STATE Update $x^{(t+1)}= x^{(t)}-p^{(t)}$;
			\ENDFOR
		\end{algorithmic}
	\end{small}
\end{algorithm}

\subsection{The Subsampled Newton method}

The	Subsampled Newton method is depicted in Algorithm~\ref{alg:H_subsamp} and the approximate Hessian is constructed by sampling:
\begin{equation}
\label{eq:samp_H}
H^{(t)} = \frac{1}{|\SM|}\sum_{j\in\mathcal{S}}\nabla^2 f_j(x^{(t)}).
\end{equation}
We now give its local convergence properties in the following theorem.
\begin{theorem}\label{thm:H_subsamp}
	Let $F(x)$ satisfy the properties described in Theorem~\ref{thm:univ_frm}.  Assume Eqn.~\eqref{eq:k} and Eqn.~\eqref{eq:sigma} hold and let $0<\delta<1$, $0<\epsilon_0<1/2$ and $0\leq\epsilon_1<1$ be given. 
	The sample size $|\SM|$ satisfies $|\SM| \geq \frac{3K/\mu\log(2d/\delta)}{\epsilon_0^2}$. 
	The approximate Hessian $H^{(t)}$ is constructed as Eqn.~\eqref{eq:samp_H}, and the direction vector $p^{(t)}$ satisfies Eqn.~\eqref{eq:inexact_p}. 
	Then for $t = 1,\dots, T$, Algorithm~\ref{alg:H_subsamp} has the following convergence properties:
	\begin{enumerate}[label = (\alph*)]
		\item There exists a sufficient small value $\gamma$ and $\nu = o(1)$ such that when $\|x^{(t)} - x^*\|_M\leq \gamma$,  each iteration satisfies  Eqn.~\eqref{eq:lin_conv} with probability at least $1-\delta$.
		\item If $\nabla^{2}F(x^{(t)})$ is also Lipschitz continuous and $\{x^{(t)}\}$ satisfies Eqn.~\eqref{eq:qua_cond}, then each iteration satisfies Eqn.~\eqref{eq:lin_qua_conv} with probability at least $1-\delta$.
		\item  If $F(x)$ is furthermore self-concordant, the iteration complexity of the sketch Newton with backtracking line search (Algorithm~\ref{alg:damp_app_newton} with $H^{(t)}$ constructed as Eqn.~\eqref{eq:samp_H}) is upper bounded by Eqn.~\eqref{eq:iter_comp}.
	\end{enumerate}
\end{theorem}
As we can see, Algorithm~\ref{alg:H_subsamp} almost has the same convergence properties as Algorithm~\ref{alg:sketch_newton} except several minor differences. The main difference is  the construction manner of  $H^{(t)}$ which should satisfy Eqn.~\eqref{eq:prec_cond}. Algorithm~\ref{alg:H_subsamp} relies on the assumption that each $\|\nabla^2f_i(x)\|$ is upper bounded (i.e., Eqn.~\eqref{eq:k} holds), while Algorithm~\ref{alg:sketch_newton} is built on the setting of the Hessian matrix as in Eqn.~\eqref{eq:H_ass}.  

\subsection{Regularized Subsampled Newton} \label{subsec:resub}

In ill-conditioned cases (i.e., $\hat{\kappa} = \frac{K}{\mu}$ is large), the subsampled Newton method  in Algorithm~\ref{alg:H_subsamp} should take a lot of samples  because the sample size $|\mathcal{S}|$  depends on $\frac{K}{\mu}$ linearly. 
To overcome this problem, one resorts to a regularized subsampled Newton method which  adds a regularizer to the original subsampled Hessian:
\begin{equation}
\label{eq:samp_H_reg}
H^{(t)} = \frac{1}{|\SM|}\sum_{j\in\mathcal{S}}\nabla^2 f_j(x^{(t)}) + \xi \cdot I
\end{equation}
where $\xi > 0$ is the regularization parameter. 
The detailed algorithmic procedure of the regularized subsampled Newton is described in Algorithm~\ref{alg:reg_subsamp}. 
In the following analysis, we  prove that adding a regularizer is an effective way to reduce the sample size while keeping converging in theory. 
\begin{theorem} \label{thm:Reg_subnewton}
Let $F(x)$ satisfy the properties described in Theorem~\ref{thm:univ_frm}.  
Assume Eqn.~\eqref{eq:k} and~\eqref{eq:sigma} hold, and let $0<\delta<1$, $0\leq\epsilon_1<1$	and $0<\xi$ be given. 
Assume the sample size $|\SM|$ satisfy $|\SM| \geq \frac{18K\log(2d/\delta)}{\xi}$, 
and $H^{(t)}$ is constructed as in Algorithm~\ref{alg:reg_subsamp}.  
	Define 
	\begin{equation}
	\epsilon_0 = \max\left(\frac{3\xi + \mu}{3\xi + 3\mu} ,\frac{L - 2\xi}{2(L+\xi)}\right),  \label{eq:eps_reg}
	\end{equation}
	which implies that $0<\epsilon_0<1$.
	Moreover, the direction vector $p^{(t)}$ satisfies Eqn.~\eqref{eq:inexact_p}. Then Algorithm~\ref{alg:reg_subsamp} has the following convergence properties:
	\begin{enumerate}[label = (\alph*)]
		\item There exists a sufficient small value $\gamma$ and $\nu = o(1)$ such that when $\|x^{(t)} - x^*\|_M\leq \gamma$,  each iteration satisfies  Eqn.~\eqref{eq:lin_conv} with probability at least $1-\delta$.
		\item If $\nabla^{2}F(x^{(t)})$ is also Lipschitz continuous and $\{x^{(t)}\}$ satisfies Eqn.~\eqref{eq:qua_cond}, then each iteration satisfies Eqn.~\eqref{eq:lin_qua_conv} with probability at least $1-\delta$.
		\item  If $F(x)$ is furthermore self-concordant, the iteration complexity of the sketch Newton with backtracking line search (Algorithm~\ref{alg:damp_app_newton} with $H^{(t)}$ constructed as Eqn.~\eqref{eq:samp_H_reg}) is upper bounded by Eqn.~\eqref{eq:iter_comp}.
	\end{enumerate}
\end{theorem}

In Theorem~\ref{thm:Reg_subnewton}	the parameter $\epsilon_0$  mainly decides  convergence properties of Algorithm~\ref{alg:reg_subsamp}. It is determined by two terms just as shown in Eqn.~\eqref{eq:eps_reg}. These two terms depict the relationship among the sample size, regularizer $\xi\cdot I$, and convergence rate. 

We can observe that the sample size $|\SM| = \frac{18K\log(2d/\delta)}{\xi}$ decreases as $\xi$ increases. 
Hence Theorem~\ref{thm:Reg_subnewton} gives a theoretical guarantee that adding the regularizer $\xi \cdot I$ is an effective approach for reducing the sample size when $K/\mu$ is large. 
Conversely, if we want to sample a small part of $f_i$'s, then we should choose a large $\xi$. 

Though a large $\xi$ can reduce the sample size, it is at the expense of a slower convergence rate. 
As we can see, $\frac{3\xi + \mu}{3\xi + 3\mu}$ goes to $1$ as $\xi$ increases. 
At the same time, $\epsilon_1$ also has to decrease. Otherwise, $\epsilon_0+\epsilon_1$ may be beyond $1$ which means that Algorithm~\ref{alg:reg_subsamp} will not converge.

In fact, a slower convergence rate in regularized subsampled Newton method is because the sample size becomes small, which implies less curvature information is obtained. 
However, a small sample size implies  low computational cost in each iteration. Therefore, a proper regularizer which balances the cost of each iteration and convergence rate is the key in the regularized subsampled Newton algorithm.  

\begin{algorithm}[tb]
\caption{Regularized Subsample Newton.}
\label{alg:reg_subsamp}
\begin{small}
	\begin{algorithmic}[1]
		\STATE {\bf Input:} $x^{(0)}$, $0<\delta<1$, regularizer parameter $\alpha$, sample size $|\SM|$ ;
		\FOR {$t=0,1,\dots$ until termination}
		\STATE Select a sample set $\SM$ of size $|\SM|$, and $H^{(t)} = \frac{1}{|\SM|}\sum_{j\in\mathcal{S}}\nabla^2 f_j(x^{(t)}) + \alpha I$;
		\STATE Calculate $p^{(t)} \approx \argmin_{p} \frac{1}{2} p^T H^{(t)} p - p^T\nabla F(x^{(t)})$	
		\STATE Update $x^{(t+1)}= x^{(t)}-p^{(t)}$;
		\ENDFOR
	\end{algorithmic}
\end{small}
\end{algorithm}

\subsection{NewSamp}

\citet{erdogdu2015convergence} proposed NewSamp which is another regularized subsampled Newton method. 
NewSamp constructs its approximate Hessian as follows:
\begin{equation}
\label{eq:newsamp}
H^{(t)} = H_{|\SM|}^{(t)} + U_{{\setminus}r}( 
\hat{\lambda}_{r+1}^{(t)}I-\hat{\Lambda}_{{\setminus}r})U_{{\setminus}r}^T ,
\end{equation}
where 
\[
H_{|\SM|}^{(t)} = \frac{1}{|\SM|}\sum_{j\in\mathcal{S}}\nabla^2 f_j(x^{(t)}), 
\]
and its SVD decomposition is 
\begin{equation*}
H_{|\SM|}^{(t)} = U\hat{\Lambda}U^T = U_r \hat{\Lambda}_r U^T_r + U_{{\setminus}r} \hat{\Lambda}_{{\setminus}r}U_{{\setminus}r}^T.
\end{equation*}
The detailed algorithm is depicted in Algorithm~\ref{alg:NewSamp}.

Now, we give the theoretical analysis of local convergence properties of NewSamp (Algorithm~\ref{alg:NewSamp}).

\begin{theorem} \label{thm:NewSamp}
Let $F(x)$ satisfy the properties described in Theorem~\ref{thm:univ_frm}.  Assume Eqn.~\eqref{eq:k} and Eqn.~\eqref{eq:sigma} hold and let $0<\delta<1$ and target rank $r$ be given. 
Let $\lambda_{r+1}$ be the $(r+1)$-th eigenvalue of $\nabla^2F(x^{(t)})$. Set the sample size $|\SM| \geq \frac{18K\log(2d/\delta)}{\lambda_{r+1}}$, and define 
\begin{equation}
\epsilon_0 = \max\left( \frac{5\lambda_{r+1} + \mu}{5\lambda_{r+1} + 3\mu},\frac{1}{2}\right), \label{eq:eps_newsamp}
\end{equation}
which implies $0<\epsilon_0<1$.
Assume the direction vector $p^{(t)}$ satisfies Eqn.~\eqref{eq:inexact_p}. Then for $t = 1,\dots, T$, Algorithm~\ref{alg:NewSamp} has the following convergence properties:
\begin{enumerate}[label = (\alph*)]
	\item There exists a sufficient small value $\gamma$ and $\nu = o(1)$ such that when $\|x^{(t)} - x^*\|_M\leq \gamma$,  each iteration satisfies  Eqn.~\eqref{eq:lin_conv} with probability at least $1-\delta$.
	\item If $\nabla^{2}F(x^{(t)})$ is also Lipschitz continuous and $\{x^{(t)}\}$ satisfies Eqn.~\eqref{eq:qua_cond}, then each iteration satisfies Eqn.~\eqref{eq:lin_qua_conv} with probability at least $1-\delta$.
	\item  If $F(x)$ is furthermore self-concordant, the iteration complexity of the sketch Newton with backtracking line search (Algorithm~\ref{alg:damp_app_newton} with $H^{(t)}$ constructed as Eqn.~\eqref{eq:newsamp}) is upper bounded by Eqn.~\eqref{eq:iter_comp}.
\end{enumerate}
\end{theorem}

\begin{algorithm}[t]
	\caption{NewSamp.}
	\label{alg:NewSamp}
	\begin{small}
		\begin{algorithmic}[1]
			\STATE {\bf Input:} $x^{(0)}$, $0<\delta<1$, $r$, sample size $|\SM|$;
			\FOR {$t=0,1,\dots$ until termination }
			\STATE Select a sample set $\SM$ of size $|\SM|$, and get $H_{|\SM|}^{(t)} = \frac{1}{|\SM|}\sum_{j\in\mathcal{S}}\nabla^2 f_j(x^{(t)})$;
			\STATE Compute rank $r+1$ truncated SVD deompostion of $H_{|\SM|}^{(t)}$ to get $U_{r+1}$ and $\hat{\Lambda}_{r+1}$. Construct $H^{(t)} = H_{|\SM|}^{(t)} + U_{{\setminus}r}( 
			\hat{\lambda}_{r+1}^{(t)}I-\hat{\Lambda}_{{\setminus}r})U_{{\setminus}r}^T$
			\STATE Calculate $p^{(t)} \approx \argmin_{p} \frac{1}{2} p^T H^{(t)} p - p^T\nabla F(x^{(t)})$	
			\STATE Update $x^{(t+1)}= x^{(t)}-p^{(t)}$;
			\ENDFOR
		\end{algorithmic}
	\end{small}
\end{algorithm}	

The first term of right hand of Eqn.~\eqref{eq:eps_newsamp} reveals the the relationship between the target rank $r$ and sample size.  
We can observe the sample size is linear to $1/\lambda_{r+1}$. 
Hence, a small $r$ means that a small sample size is sufficient. 
Conversely, if we want to sample a small portion of $f_i$'s, then we should choose a small $r$. 
Eqn.~\eqref{eq:eps_newsamp} shows that a small sample size will lead to a poor convergence rate. 
If we set $r = 0$, then $\epsilon_0$ will be $1-\frac{2\mu}{5\lambda_1+3\mu}$. Consequently, the convergence rate of NewSamp is almost the same as gradient descent. 

It is worth pointing out that Theorem~\ref{thm:NewSamp} explains the empirical results that NewSamp is applicable in training SVM in which the Lipschitz continuity condition of $\nabla^2F(x)$ is not satisfied \cite{erdogdu2015convergence}. 

\begin{table*}
	\centering
	\caption{Comparison with previous work. We use (Reg)SSN to denote the (regularized) subsampled Newton method. For ReSNN, we choose $\xi = \lambda_{r+1}$. The notation $\ti{\OM}(\cdot)$ hides the polynomial of $\log (d/\delta)$.}
	\label{tb:comp_ssn}
	\begin{tabular}{c|ccc}
		\hline
		~~~~Method~~~~& ~~~~Reference ~~~~& ~~~~Sample Size ~~~~& ~~~~Iterations Complexity~~~~ \\ \hline
		\multirow{2}*{SSN} 
		&Theorem 5 of \cite{roosta2019sub}			& $\tilde{\OM}(K^2/\mu^2)$	& $\tilde{\OM}(\log(1/\epsilon))$ \\
		&{{\bf Theorem~\ref{thm:H_subsamp}}} 	&$\tilde{\OM}(K/\mu)$			& $\tilde{\OM}(\log(1/\epsilon))$ \\
		\hline
		\multirow{1}*{RegSNN} 
		&{{\bf Theorem~\ref{thm:Reg_subnewton}}} & $\tilde{\OM}\left(K/\lambda_{r+1}\right)$ & $\tilde{\OM}\left(\frac{\lambda_{r+1}}{\mu}\log(1/\epsilon)\right)$\\
		\hline
		\multirow{2}*{NewSamp} 
		&Theorem 3.2 of \cite{erdogdu2015convergence} & $\tilde{\OM}\left(K^2/\mu^2\right)$  & $\OM(\log(1/\epsilon))$  \\
		&{{\bf Theorem~\ref{thm:NewSamp}}} & $\tilde{\OM}\left(K/\lambda_{r+1}\right)$ & $\OM\left(\frac{\lambda_{r+1}}{\mu}\log(1/\epsilon)\right)$\\
		\hline
	\end{tabular}
\end{table*}

\subsection{Comparison with Previous Work}

We will compare our results in this section with previous work. 
Though many variants of subsampled Newton methods have been proposed recently, 
they share the similar proof procedure. 
Thus, these algorithms have almost the same sample size and convergence rate.
For example, the subsampled Newton method \cite{roosta2019sub} and NewSamp \cite{erdogdu2015convergence} have the same order of sample size and convergence rate (referring to Table~\ref{tb:comp_ssn}).   
Thus, we only compare our results with the recent work of \citet{roosta2019sub} and NewSamp \cite{erdogdu2015convergence}.
The detailed comparison is listed in Table~\ref{tb:comp_ssn}. 

First, compare our analysis of subsampled Newton with the one of \cite{roosta2019sub}.
We can observe that to achieve the same convergence rate, our result only needs $\ti{\OM}(K/\mu)$ in contrast to $\ti{\OM}(K^2/\mu^2)$ of \citet{roosta2019sub}.
Hence, our result is substantially much tighter than previous work.

Then we compare our theoretical analysis of NewSamp with the one of \citet{erdogdu2015convergence}.
We can observe that though  NewSamp is a kind of regularized subsampled Newton, it still takes $\ti{\OM}(K^2/\mu^2)$ samples which is the same to subsampled Newton.
In contrast, our analysis (Theorem~\ref{thm:NewSamp}) describes how to the regularization reduces the sample size and convergence speed.
This theory matches the empirical study that a small $r$ (implying a large $\lambda_{r+1}$) will reduce the samples and convergence speed \cite{erdogdu2015convergence}.

Finally, we compare NewSamp with regularized subsampled Newton (Algorithm~\ref{alg:reg_subsamp}).    
We mainly focus on the parameter $\epsilon_0$ in Theorem~\ref{thm:Reg_subnewton} and Theorem~\ref{thm:NewSamp} which mainly determines convergence properties of Algorithm~\ref{alg:reg_subsamp} and Algorithm~\ref{alg:NewSamp}. 
Specifically, if we set $\xi = \lambda_{r+1}$ in Eqn.~\eqref{eq:eps_reg}, then $\epsilon_0 = \frac{3\lambda_{r+1} + \mu}{3\lambda_{r+1} + 3\mu} $ which is of the same order of the first term of the right hand of Eqn.~\eqref{eq:eps_newsamp}. Hence, we can regard NewSamp as a special case of Algorithm~\ref{alg:reg_subsamp}. However, NewSamp provides an approach  for automatic choice of $\alpha$. 
Recall that NewSamp includes another parameter: the target rank $r$. Thus, NewSamp and Algorithm~\ref{alg:reg_subsamp} have the same number of free parameters. If $r$ is not properly chosen, NewSamp will still have poor performance. Therefore, Algorithm~\ref{alg:reg_subsamp} is theoretically preferred because NewSamp needs extra cost to perform SVDs.

\section{Empirical Analysis}
\label{sec:exp}

In this section, we  validate our theoretical results about unnecessity of the Lipschitz continuity condition of $\nabla^2F(x)$, sketched size of the sketch Newton and how the regularization affects the sample size and convergence rate of regularized Newton, experimentally.. 

\begin{table}[!ht]
	\centering
	\caption{Datasets Description}
	\label{tb:data}
	\begin{tabular}{cccc}
		\hline
		Dataset & $n$        & $d$         & source \\ \hline
		mushrooms& $8,124$    & $112$       & UCI       \\  \hline 
		a9a     &  $32,561$   & $123$       & UCI   \\ \hline
		Covertype   &  $581,012$   & $54$   & UCI \\ \hline
	\end{tabular}
\end{table}

\subsection{Unnecessity of Lipschitz continuity of Hessian}

We conduct experiment on the primal problem for the linear SVM which can be written as
\begin{align*}
\min_{x}F(x) = \frac{1}{2}\|x\|^2 + \frac{C}{2n} \sum_{i=1}^{n}\ell(b_i,\langle x, a_i\rangle)
\end{align*}
where $(a_i,b_i)$ denotes the training data, $x$ defines the separating hyperplane, $C>0$, and $\ell(\cdot)$ is the loss function. In our experiment, we choose Hinge-2 loss as our loss function whose definition is 
\begin{align*}
\ell(b,\langle x, a \rangle) = \max(0, 1-b\langle x, a \rangle)^2.
\end{align*}

Let $SV^{(t)}$ denote the set of indices of all the support vectors at iteration $t$, i.e.,
\begin{align*}
SV^{(t)} = \{i: b_i\langle x^{(t)}, a_i\rangle <1 \}.
\end{align*}
Then the Hessian matrix of $F(x^{(t)})$ can be written as 
\begin{align*}
\nabla^2F(x^{(t)}) = I + \frac{1}{n}\sum_{i\in SV^{(t)}} a_ia_i^T.
\end{align*} 
From the above equation, we can see that $\nabla^2F(x)$ is not Lipschitz continuous.

\begin{figure}[!ht]
	\subfigtopskip = 0pt
	\begin{center}
		\centering
		\subfigure[\textsf{mushrooms.}]{\includegraphics[width=51mm]{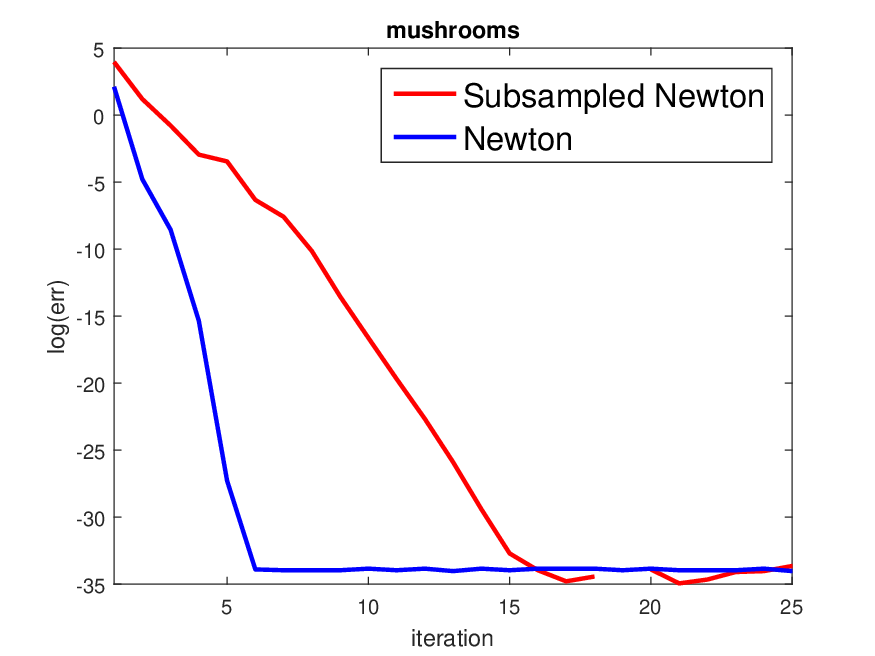}}~
		\subfigure[\textsf{a9a.}]{\includegraphics[width=51mm]{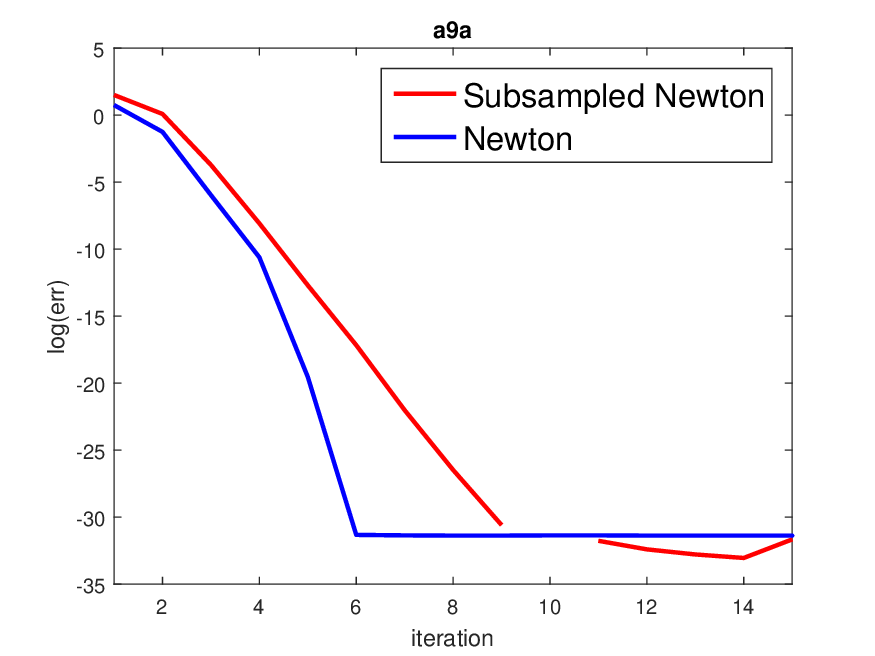}}~
		\subfigure[\textsf{covtype.}]{\includegraphics[width=51mm]{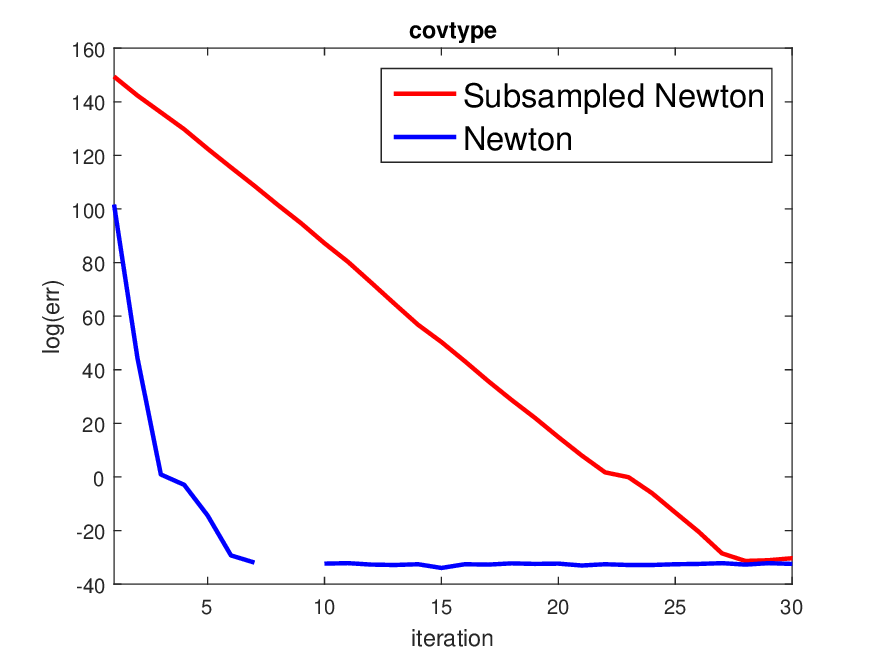}}\\
	\end{center}
	\caption{Convergence properties on different datasets.}
	\label{fig:conv}
\end{figure}

Without loss of generality, we use the Subsampled Newton method (Algorithm~\ref{alg:H_subsamp}) in our experiment. We sample $5\%$ support vectors in each iteration.  
Our experiments on three datasets whose detailed description is in Table~\ref{tb:data} and report our results in Figure~\ref{fig:conv}. 

From Figure~\ref{fig:conv}, we can see that Subsampled Newton converges linearly and the Newton method converges superlinearly. This matches our theory that the Lipschitz continuity of $\nabla^2F(x)$ is not necessary to achieve a linear or superlinear convergence rate.  

\subsection{Sketched Size of Sketch Newton}

\begin{figure}[!ht]
	\subfigtopskip = 0pt
	\begin{center}
		\centering
		\subfigure[\textsf{Gaussian Sketching.}]{\includegraphics[width=51mm]{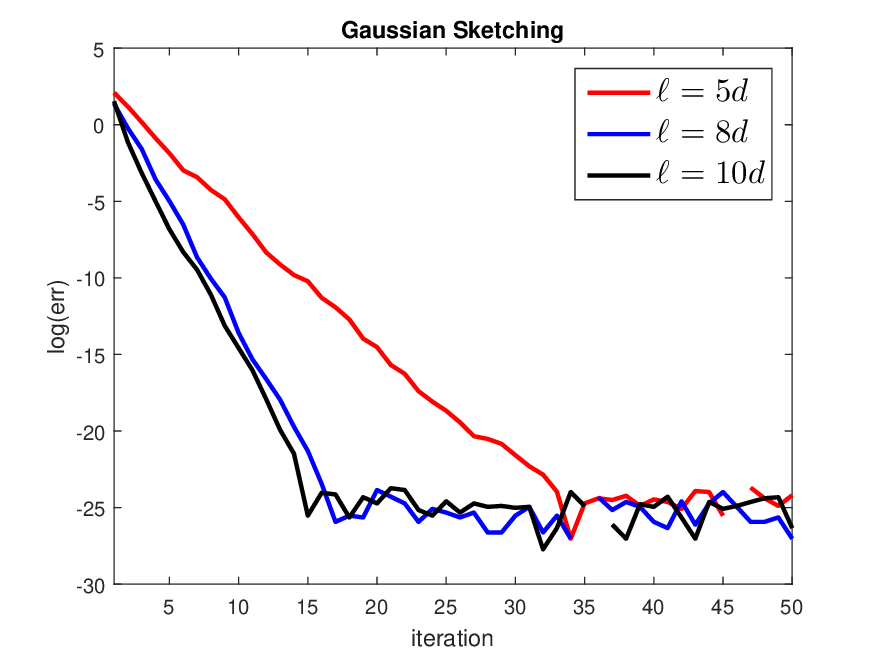}}~
		\subfigure[\textsf{Sparse Sketching.}]{\includegraphics[width=51mm]{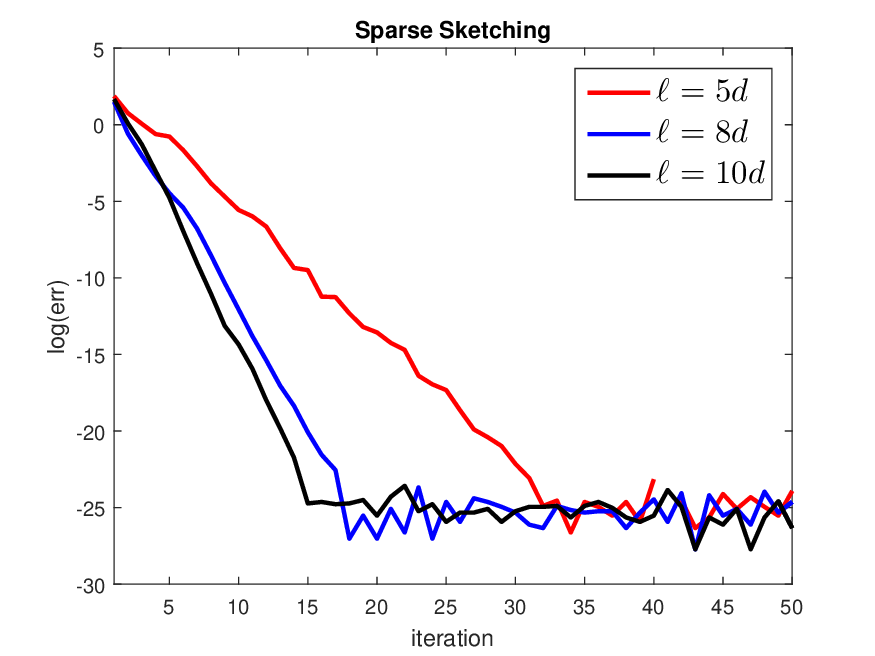}}~
		\subfigure[\textsf{Leverage Score Sketching.}]{\includegraphics[width=51mm]{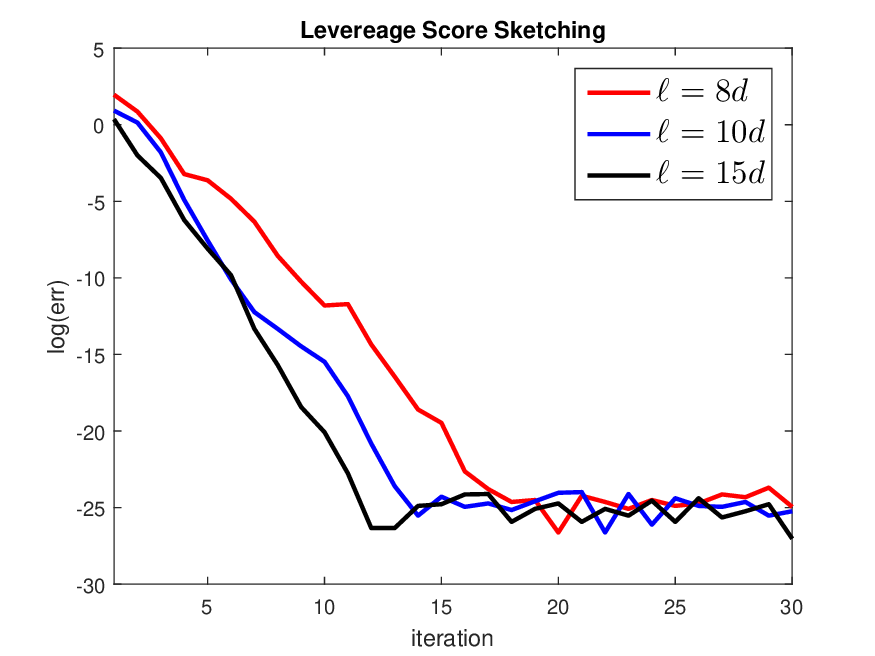}}\\
	\end{center}
	\caption{Convergence properties of different sketched sizes}
	\label{fig:sketch}
\end{figure}

Now we validate that our theoretical result that sketched size is independent of the condition number of the Hessian matrix in Sketch Newton. To control the condition number of the Hessian conveniently, we conduct the experiment on least squares regression which is defined as
\begin{align}
\min_x\frac{1}{2}\|Ax - b\|^2. \label{eq:lsr}
\end{align}
In each iteration, the Hessian matrix is $A^TA$. In our experiment, $A$ is a $10000\times 54$ matrix. And we set the singular values $\sigma_i$ of $A$ as:
\[
\sigma_i = 1.2^{-i}.
\]
Then the condition number of $A$ is $\kappa(A) = 1.2^{54} = 1.8741\times 10^{4}$. We use different sketch matrices in Sketch Newton (Algorithm~\ref{alg:sketch_newton}) and set different values of the sketched size $\ell$. We report our empirical results in Figure~\ref{fig:sketch}.

From Figure~\ref{fig:sketch}, we can see that Sketch Newton performs well when the sketch size $\ell$ is several times of $d$ for all different sketching matrices. Moreover, the corresponding algorithms converge linearly.  This matches our theory that sketched size is independent of the condition number of Hessian matrix to achieve a linear convergence rate. In contrast, the theoretical result  of \cite{xu2016sub} shows that sketched size is $\ell = d * \kappa(A) = 1.02\times 10^6$  bigger than $n= 10^4$.  

\subsection{Sample Size of Regularized Subsampled Newton}

\begin{figure}[!ht]
	\subfigtopskip = 0pt
	\begin{center}
		\centering
		\subfigure[\textsf{Sample Size $|\SM| = 100$.}]{\includegraphics[width=51mm]{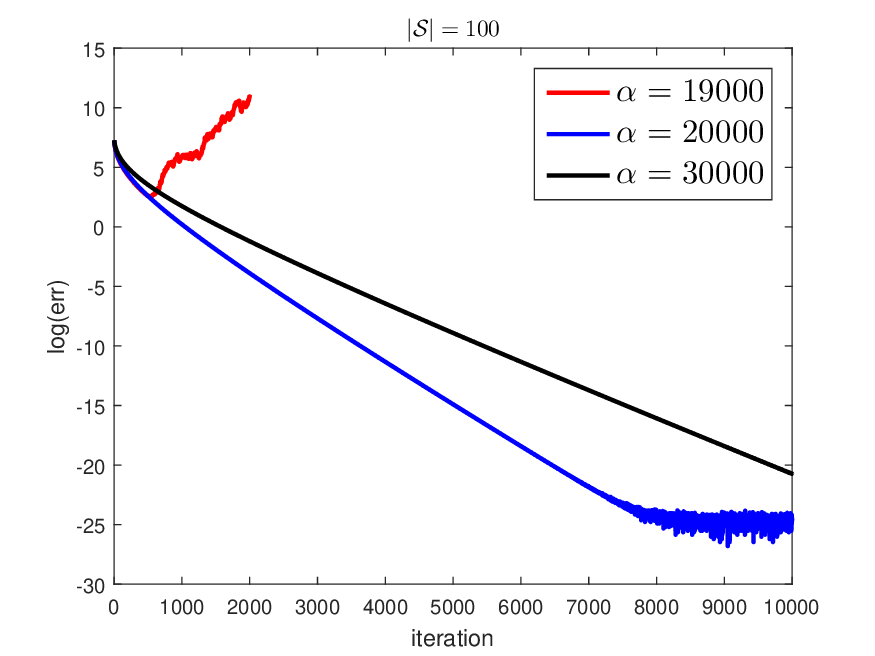}}~
		\subfigure[\textsf{Sample size $|\SM| = 300$.}]{\includegraphics[width=51mm]{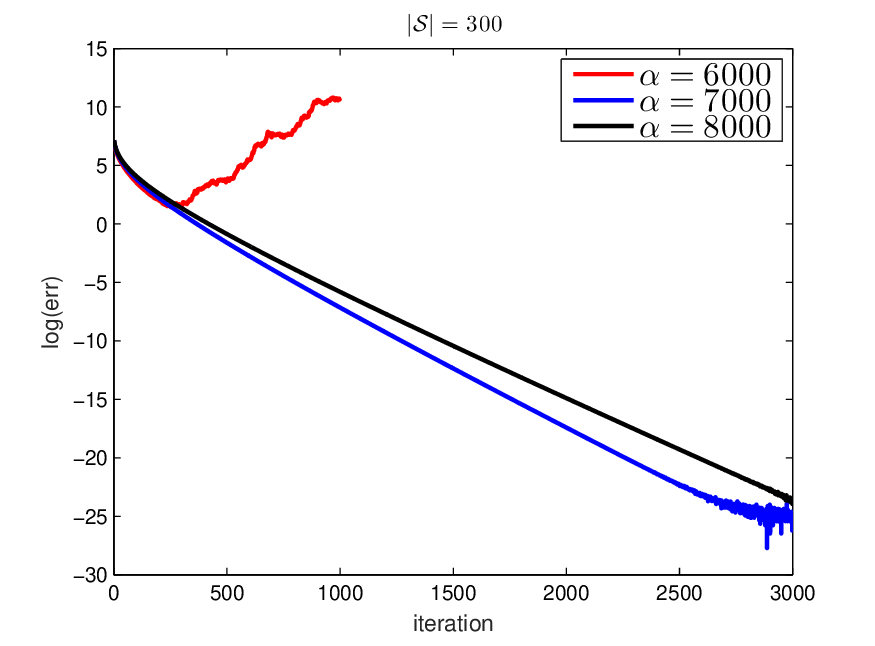}}~
		\subfigure[\textsf{Sample size $|\SM| = 600$.}]{\includegraphics[width=51mm]{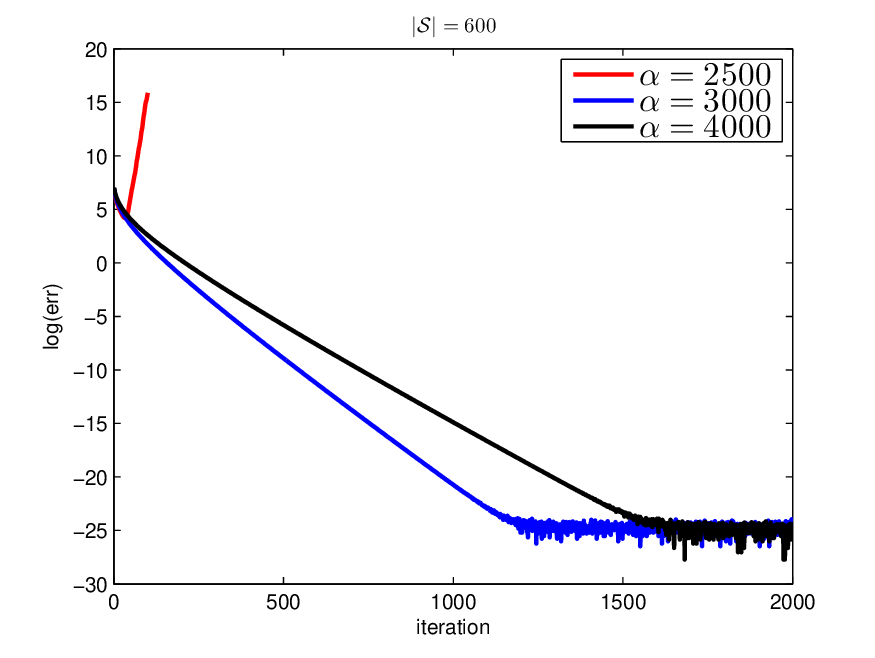}}\\
	\end{center}
	\caption{Convergence properties of Regularized Subsampled Newton}
	\label{fig:resub}
\end{figure}

\begin{figure}[!ht]
	\subfigtopskip = 0pt
	\begin{center}
		\centering
		\subfigure[\textsf{Sample Size $|\SM| = 100$.}]{\includegraphics[width=51mm]{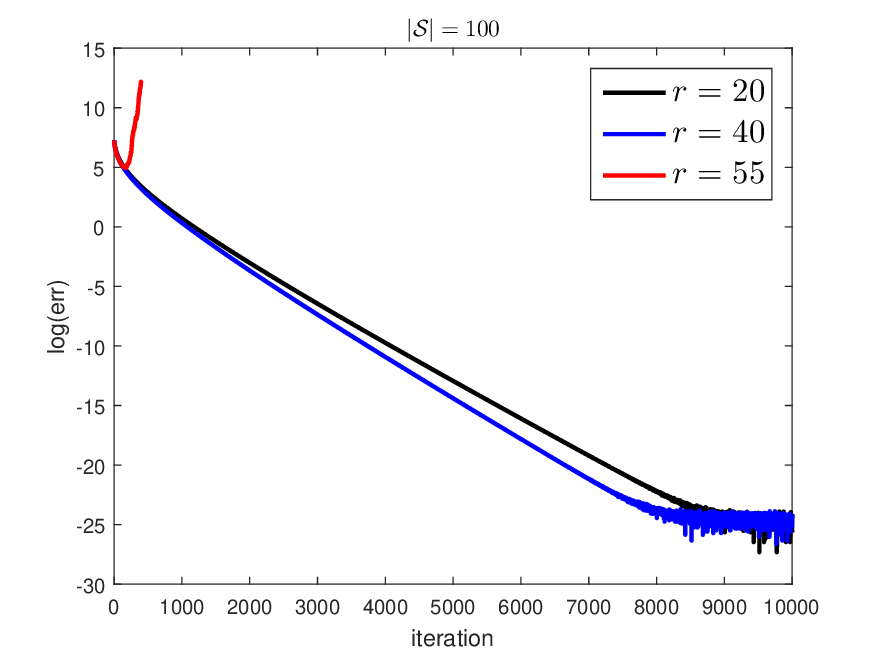}}~
		\subfigure[\textsf{Sample size $|\SM| = 300$.}]{\includegraphics[width=51mm]{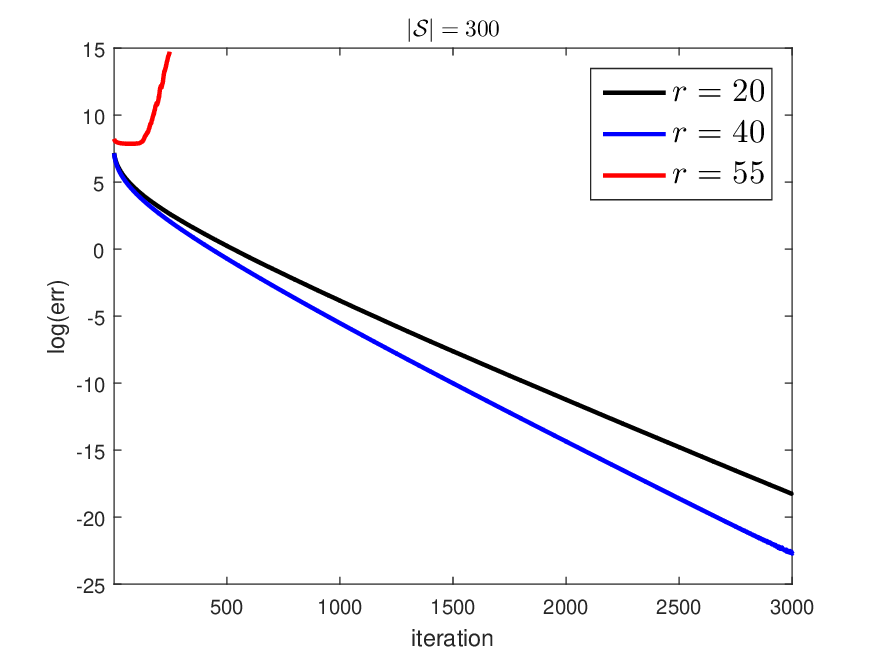}}~
		\subfigure[\textsf{Sample size $|\SM| = 600$.}]{\includegraphics[width=51mm]{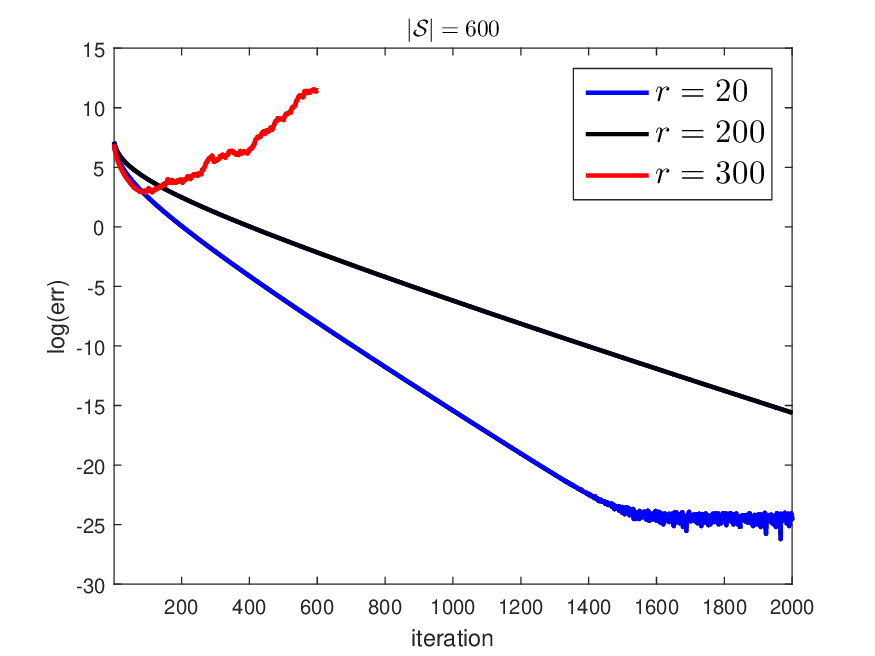}}\\
	\end{center}
	\caption{Convergence properties of NewSamp}
	\label{fig:NewSamp}
\end{figure}

We also choose least squares regression defined in Eqn.~\eqref{eq:lsr} in our experiment to validate the theory that adding a regularizer is an effective approach to reducing the sample size while keeping convergence in Subsampled Newton. Let $A\in\RB^{n\times d}$ where $n = 8000$ and $d = 5000$. Hence Sketch Newton can not be used in this case because $n$ and $d$ are close to each other. In our experiment, we set different sample sizes $|\SM|$. For each $|\SM|$ we choose different regularizer terms $\alpha$ and different target ranks $r$. We report our results in Figures~\ref{fig:resub} and~\ref{fig:NewSamp}.

As we can see, if the sample size $|\SM|$ is small, then we should choose a large $\alpha$ in Algorithm~\ref{alg:reg_subsamp}; otherwise, the algorithm will diverge. However, if the regularizer term $\alpha$ is too large, then the algorithm will converge slowly. Besides, increasing the sample size and choosing a proper regularizer will improve convergence properties obviously. When $|\SM| = 600$, it only needs about $1200$ iterations to obtain a precise solution while it needs about $8000$ iterations when $|\SM| = 100$. Similarly, if the sample size $|\SM|$ is small, then we should choose a small target rank in NewSamp. Otherwise NewSamp may diverge. Also if the target rank is not chosen properly, then NewSamp will have poor convergence properties. Furthermore, comparing Figures~\ref{fig:resub} and~\ref{fig:NewSamp}, we can see that the two algorithms have similar convergence properties. This validates the theoretical result that NewSamp provides a method to choose $\alpha$ automatically. Our empirical analysis matches the theoretical analysis in Subsection~\ref{subsec:resub} very well.    

\section{Conclusion}
\label{sec:conclusion}	

In this paper, we have proposed a framework to analyze both local and global convergence properties of second order methods including stochastic  and deterministic versions. This framework reveals some important convergence properties of the subsampled Newton method and sketch Newton method,  which are unknown before. The most important thing is in that our analysis lays the theoretical foundation of several important stochastic second order methods. 	

We believe that this framework might also provide some useful insights for developing new subsampled Newton-type algorithms.  We would like to address this issue in future.   

\bibliographystyle{plainnat}
\bibliography{ref}

\newpage
\appendix

\section{Some Important Lemmas}	In this section, we give several important lemmas which will be used in the proof of the theorems of this paper.

\begin{lemma}\label{lem:prec}
	If $A$ and $B$ are $d\times d$ symmetric positive matrices, and $(1-\epsilon_0)B \preceq A \preceq (1+\epsilon_0) B$ where $0<\epsilon_0<1$, then we have
	\[
	\|A^{1/2}B^{-1} A^{1/2} - I\| \leq \epsilon_0,
	\]
	where $I$ is the identity matrix.
\end{lemma}
\begin{proof}
	Because $A \preceq (1+\epsilon_0) B$, we have $z^T [A- (1+\epsilon_0) B] z \leq 0$ for any nonzero $z \in \RB^d$. This implies $ \frac{z^T A z}{z^T B z} \leq 1+ \epsilon_0$ for any $z \neq 0$. Subsequently, 
	\begin{align*}
	\lambda_{\max}(B^{-1}A) =& \lambda_{\max}(B^{-1/2}AB^{-1/2})\\
	=& \max_{u \neq 0} \frac{u^T B^{-1/2}AB^{-1/2}u}{u^T u}\\
	=& \max_{z \neq 0} \frac{z^T A z}{z^T B z}\\
	\leq& 1+\epsilon_0, 
	\end{align*}
	where the last equality is obtained by setting $z = B^{-1/2}u$. Similarly, we have $\lambda_{\min}(B^{-1}A) \geq 1-\epsilon_0$.	
	Since $B^{-1}A$ and $A^{1/2}B^{-1} A^{1/2}$ are similar, the eigenvalues of $A^{1/2}B^{-1} A^{1/2}$ are all between $1-\epsilon_0$ and $1+\epsilon_0$. Therefore, we have 
	\[
	\|A^{1/2}B^{-1} A^{1/2} - I\| \leq \epsilon_0.
	\]	
\end{proof}

\begin{lemma}[\cite{tropp2015introduction}] \label{lem:matrix_bnd}
	Let $X_{1},X_2,\dots, X_k$ be independent, random, symmetric, real matrices of size $d\times d$ with $0\preceq X_i \preceq L I$, where $I$ is the $d\times d$ identity matrix. Let $Y = \sum_{i=1}^{k}{X_i}$, $\mu_{\min} = \lambda_{\min}(\EB[Y])$ and $\mu_{\max} = \lambda_{\max}(\EB[Y])$. Then, we have
	\[
	\PB\left(\lambda_{\min}(Y)\leq (1-\epsilon)\mu_{\min}\right) \leq d\cdot e^{-\epsilon^{2}\mu_{\min}/2L},
	\]
	and
	\[
	\PB\left(\lambda_{\max}(Y)\geq (1+\epsilon)\mu_{\max}\right) \leq d\cdot e^{-\epsilon^{2}\mu_{\min}/3L}.
	\]
\end{lemma}

\section{Proofs of Theorem~\ref{thm:univ_frm}} \label{sec:app_proof}

The proof Theorem~\ref{thm:univ_frm} consists of the following lemmas. 
First, by Lemma~\ref{lem:decr}, we upper bound  $\norm{x^{(t+1)} - x^*}_M$ by three terms.
The first term dominates the convergence property. 
The second term depicts how the approximate descent direction affects the convergence.
The third term is a high order term.
 
In Lemma~\ref{lem:cond}, we prove that the first term of right hand of Eqn.~\eqref{eq:upper} is upper bounded by $\epsilon_0 \norm{x^{(t)} - x^*}$ and a high order term.
Lemma~\ref{lem:inexact} shows that the second term affect the convergence rate at most $\epsilon_1$.
In Lemma~\ref{lem:un_lip}, we complete the convergence analysis when the Hessian is continuous near the optimal point but the Hessian is not Lipschitz continuous.
If the the Hessian is not Lipschitz continuous, Lemma~\ref{lem:lip} provides the detailed convergence analysis.

\begin{lemma}
	\label{lem:decr}
	Letting sequence $\{x^{(t)}\}$ update as Algorithm~\ref{alg:app_newton}, then it satisfies 
	\begin{equation}
	\label{eq:upper}
	\small
	\begin{aligned}
	\norm{x^{(t+1)} - x^*}_M \leq& \norm{I - M^{1/2}[H^{(t)}]^{-1}M^{1/2}}\cdot\norm{x^{(t)} - x^*}_M + \norm{[H^{(t)}]^{-1}\nabla F(x^{(t)}) - p^{(t)}}_M\\
	&+\norm{M^{1/2}[H^{(t)}]^{-1} \left(\nabla F(x^{(t)}) - \nabla^2 F(x^*) (x^{(t)} - x^*) \right) }
	\end{aligned}
	\end{equation}
	where $M = \nabla^2 F(x^*)$.
\end{lemma}
\begin{proof}
	By the update procedure of $x^{(t)}$, we have
	\begin{align*}
	x^{(t+1)} - x^* =& x^{(t)} - x^* - p^{(t)}\\
	=& x^{(t)} - x^* - [H^{(t)}]^{-1}\nabla F(x^{(t)}) +  [H^{(t)}]^{-1}\nabla F(x^{(t)}) - p^{(t)} \\
	=&x^{(t)} - x^* +  [H^{(t)}]^{-1}\nabla F(x^{(t)}) - p^{(t)} \\& - [H^{(t)}]^{-1} \left(\nabla^2 F(x^*) (x^{(t)} - x^*) + \nabla F(x^{(t)}) - \nabla^2 F(x^*) (x^{(t)} - x^*)  \right).
	\end{align*}
	Letting us denote $M = \nabla^2 F(x^*)$, and multiplying $M^{1/2}$ to the left and right hands of above equality, we can obtain that
	\begin{align*}
	M^{1/2} (x^{(t+1)} - x^*) 
	=& M^{1/2} (x^{(t)} - x^*) - M^{1/2}[H^{(t)}]^{-1}M^{1/2}\cdot M^{1/2} (x^{(t)} - x^*)\\
	&+M^{1/2} \left( [H^{(t)}]^{-1}\nabla F(x^{(t)}) - p^{(t)}\right)\\
	&-M^{1/2}[H^{(t)}]^{-1} \left(\nabla F(x^{(t)}) - \nabla^2 F(x^*) (x^{(t)} - x^*) \right).
	\end{align*}
	Thus, we can obtain that 
	\begin{align*}
	\norm{x^{(t+1)} - x^*}_M \leq& \norm{I - M^{1/2}[H^{(t)}]^{-1}M^{1/2}}\cdot\norm{x^{(t)} - x^*}_M + \norm{[H^{(t)}]^{-1}\nabla F(x^{(t)}) - p^{(t)}}_M\\
	&+\norm{M^{1/2}[H^{(t)}]^{-1} \left(\nabla F(x^{(t)}) - \nabla^2 F(x^*) (x^{(t)} - x^*) \right) }
	\end{align*}
\end{proof}

\begin{lemma}
	\label{lem:cond}
	Assume that the objective function $F(x)$ satisfies Assumption 1 and 2. Let $M$ denote $\nabla^2F(x^*)$, and the approximate Hessian $H^{(t)}$ satisfy Condition~\eqref{eq:prec_cond}. Then if  $\norm{\Delta}$ is sufficient small with
	\begin{equation}
	\Delta = \nabla^2F(x^*) - \nabla^2F(x^{(t)}),
	\end{equation} 
	we have 
	\begin{equation*}
	\norm{I- M^{1/2}[H^{(t)}]^{-1}M^{1/2}} \leq \epsilon_0+\mu^{-1/2}\norm{\Delta}^{1/2}(1+\epsilon_0)\left(2+\mu^{-1/2}\norm{\Delta}^{1/2}\right).
	\end{equation*}
\end{lemma}
\begin{proof}
	If $\norm{\Delta}$ is sufficient small (which implies that $\nabla^2F(x^*)$ and $\nabla^2F(x^{(t)})$ are close enough), then we have
	\begin{align*}
	\lambda_{\max} \left([\nabla^2F(x^*)]^{\frac{1}{2}}[H^{(t)}]^{-1}[\nabla^2F(x^*)]^{\frac{1}{2}}\right) = 1+\epsilon_0'\\
	\lambda_{\min} \left([\nabla^2F(x^*)]^{\frac{1}{2}}[H^{(t)}]^{-1}[\nabla^2F(x^*)]^{\frac{1}{2}}\right) = 1 - \epsilon_0''
	\end{align*}
	with $0< \epsilon_0'<1,0<\epsilon_0''<1$.
	
	Now we consider the case 
	\begin{equation}
	\label{eq:pos}
	\norm{I- [\nabla^2F(x^*)]^{\frac{1}{2}}[H^{(t)}]^{-1}[\nabla^2F(x^*)]^{\frac{1}{2}}} = 
	\epsilon_0'
	\end{equation} 
	which implies $\epsilon_0'\geq \epsilon_0''$. 
	By the properties of eigenvalue and singular value of matrices, we have
	\begin{align*}
	\lambda^2_{\max} \left(M^{\frac{1}{2}}[H^{(t)}]^{-1}M^{\frac{1}{2}}\right) = \lambda^2_{\max} \left([H^{(t)}]^{-1}M\right) \leq \sigma^2_{1}\left([H^{(t)}]^{-1}M\right)= \lambda_{\max} \left(M[H^{(t)}]^{-2}M\right),
	\end{align*}
	where the inequality follows from the fact that the largest eigenvalue is no larger than the largest singular value. 
	Thus, we obtain that
	\begin{equation*}
	\epsilon_0' = \lambda_{\max} \left(M^{\frac{1}{2}}[H^{(t)}]^{-1}M^{\frac{1}{2}}\right) \leq {\lambda^{1/2}_{\max}\left(M[H^{(t)}]^{-2}M\right)}
	\end{equation*}
	Since Eqn.~\eqref{eq:pos} holds, then we have
	\begin{equation}
	\label{eq:up_3}
	\norm{I- M^{1/2}[H^{(t)}]^{-1}M^{1/2}} \leq \norm{I- \left(M[H^{(t)}]^{-2}M \right)^{1/2}}
	\end{equation}
	
	Next, we will prove that Eqn.~\eqref{eq:up_3} still holds when $\epsilon_0'<\epsilon_0''$ which will lead to 
	\begin{equation*}
	\norm{I- [\nabla^2F(x^*)]^{\frac{1}{2}}[H^{(t)}]^{-1}[\nabla^2F(x^*)]^{\frac{1}{2}}} = 
	\epsilon_0''.
	\end{equation*}
	By the properties of eigenvalue and singular value of matrices, we have
	\begin{align*}
	\lambda^2_{\min} \left(M^{\frac{1}{2}}[H^{(t)}]^{-1}M^{\frac{1}{2}}\right) = \lambda^2_{\min} \left([H^{(t)}]^{-1}M\right) \geq \sigma^2_{\min}\left([H^{(t)}]^{-1}M\right)= \lambda_{\min} \left(M[H^{(t)}]^{-2}M\right),
	\end{align*}
	where the inequality follows from the fact that the smallest eigenvalue is no smaller than the smallest singular value. 
	This implies that
	\begin{equation*}
	\epsilon_0'' = \lambda_{\min} \left(M^{\frac{1}{2}}[H^{(t)}]^{-1}M^{\frac{1}{2}}\right) \geq  {\lambda^{1/2}_{\min}\left(M[H^{(t)}]^{-2}M\right)}
	\end{equation*}
	which implies that Eqn.~\eqref{eq:up_3} holds.
	
	Next, we will upper bound the value of right hand of Eqn.~\eqref{eq:up_3}. First, we consider the case that 
	\begin{equation}
	\label{eq:cond}
	\lambda_{\max}\left(M[H^{(t)}]^{-2}M \right)^{1/2} - 1 \geq 1 - \lambda_{\min} \left(M[H^{(t)}]^{-2}M \right)^{1/2},
	\end{equation}
	which implies that 
	\begin{equation*}
	\norm{I- \left(M[H^{(t)}]^{-2}M \right)^{1/2}} = \lambda_{\max}\left(M[H^{(t)}]^{-2}M \right)^{1/2} - 1.
	\end{equation*}
	Furthermore, we have
	\begin{align*}
	\lambda_{\max}\left(M[H^{(t)}]^{-2}M \right)^{1/2} - 1=&\norm{M[H^{(t)}]^{-2}M}^{1/2} - 1\\
	=& \norm{M[H^{(t)}]^{-2}M}^{1/2} - \norm{\tI}^{1/2} + \norm{\tI}^{1/2} - 1\\
	\leq& \epsilon_0+ \norm{M[H^{(t)}]^{-2}M}^{1/2} - \norm{\tI}^{1/2} 
	\end{align*}
	where we denote 
	\begin{equation*}
		\tI = \nabla^2 F(x^{(t)}) [H^{(t)}]^{-2} \nabla^2 F(x^{(t)}),
	\end{equation*}
	and the last inequality follows the condition~\eqref{eq:prec_cond}.
	
	Moreover, we have
	\begin{align*}
	&\norm{M[H^{(t)}]^{-2}M}^{1/2} - \norm{\tI}^{1/2} \\
	=& \norm{\tI + \Delta[H^{(t)}]^{-2} \nabla^2 F(x^{(t)})
		+[H^{(t)}]^{-2} \nabla^2 F(x^{(t)})\Delta +  \Delta [H^{(t)}]^{-2} \Delta}^{1/2} - \norm{\tI}^{1/2}\\
	\leq&\norm{\tI}^{1/2} + \norm{\Delta[H^{(t)}]^{-2} \nabla^2 F(x^{(t)})
		+[H^{(t)}]^{-2} \nabla^2 F(x^{(t)})\Delta +  \Delta [H^{(t)}]^{-2} \Delta}^{1/2} - \norm{\tI}^{1/2}\\
	\leq&2\norm{\Delta}^{1/2}\norm{[H^{(t)}]^{-1/2}}\cdot \norm{[H^{(t)}]^{-1}\nabla^2 F(x^{(t)})}^{1/2} + \norm{[H^{(t)}]^{-1}}\cdot\norm{\Delta}.
	\end{align*}
	By Condition~\eqref{eq:prec_cond}, we can obtain that
	\begin{equation}
	\norm{[H^{(t)}]^{-1}} \leq (1+\epsilon_0)\norm{[\nabla^2 F(x^{(t)})]^{-1}} \leq (1+\epsilon_0)\mu^{-1}
	\end{equation}
	and
	\begin{equation}
	\norm{[H^{(t)}]^{-1}\nabla^2 F(x^{(t)})} = \lambda^{1/2}_{\max}\left(\nabla^2 F(x^{(t)}) [H^{(t)}]^{-2}\nabla^2 F(x^{(t)})\right)\leq (1+\epsilon_0).
	\end{equation}
	Thus, we can obtain that
	\begin{equation*}
	\norm{M[H^{(t)}]^{-2}M}^{1/2} - \norm{\tI}^{1/2} 
	\leq \mu^{-1/2}\norm{\Delta}^{1/2}(1+\epsilon_0)\left(2+\mu^{-1/2}\norm{\Delta}^{1/2}\right).
	\end{equation*}
	
	Now we consider the case that Eqn.~\eqref{eq:cond} does not hold which implies that
	\begin{equation*}
	\norm{I- \left(M[H^{(t)}]^{-2}M \right)^{1/2}} = 1 - \lambda_{\min}\left(M[H^{(t)}]^{-2}M \right)^{1/2}.
	\end{equation*}
	Furthermore, we have
	\begin{align*}
	1 - \lambda_{\min}\left(M[H^{(t)}]^{-2}M \right)^{1/2} =&1 - \lambda^{1/2}_{\min}(\tI) + \lambda^{1/2}_{\min}(\tI) - \lambda^{1/2}_{\min}\left(M[H^{(t)}]^{-2}M\right)\\
	\leq& \epsilon_0 + \lambda^{1/2}_{\min}(\tI) - \lambda^{1/2}_{\min}\left(M[H^{(t)}]^{-2}M\right),
	\end{align*}
	where the last inequality follows from condition~\eqref{eq:prec_cond}.
	Since $\norm{\Delta}$ is sufficient small, then we have that
	\begin{align*}
	&\lambda^{1/2}_{\min}(\tI) - \lambda^{1/2}_{\min}\left(M[H^{(t)}]^{-2}M\right)\\ 
	=&\lambda^{1/2}_{\min}(\tI) - \lambda^{1/2}_{\min}\left(\tI + \Delta[H^{(t)}]^{-2} \nabla^2 F(x^{(t)})
	+[H^{(t)}]^{-2} \nabla^2 F(x^{(t)})\Delta +  \Delta [H^{(t)}]^{-2} \Delta\right)\\
	\leq&\lambda^{1/2}_{\min}(\tI) - \lambda^{1/2}_{\min} (\tI) + \norm{\Delta[H^{(t)}]^{-2} \nabla^2 F(x^{(t)})
		+[H^{(t)}]^{-2} \nabla^2 F(x^{(t)})\Delta +  \Delta [H^{(t)}]^{-2} \Delta}^{1/2} \\
	\leq& \mu^{-1/2}\norm{\Delta}^{1/2}(1+\epsilon_0)\left(2+\mu^{-1/2}\norm{\Delta}^{1/2}\right),
	\end{align*}
	where  the first inequality is because of $\lambda_{\min}(A+B) = \sigma_{\min}(A+B) \geq \sigma_{\min}(A) - \norm{B}$ and the fact that $(a -b)^{1/2} \geq a^{1/2} - b^{1/2}$ if $a\geq b$ and $a,b\geq 0$.
	
	Therefore, we can obtain that
	\begin{equation*}
	\norm{I- M^{1/2}[H^{(t)}]^{-1}M^{1/2}} \leq \epsilon_0+\mu^{-1/2}\norm{\Delta}^{1/2}(1+\epsilon_0)\left(2+\mu^{-1/2}\norm{\Delta}^{1/2}\right).
	\end{equation*}
\end{proof}

\begin{lemma}
	\label{lem:inexact}
	Let $p^{(t)}$ satisfy Condition~\eqref{eq:inexact_p} and $F(x)$ satisfy Assumption 1 and 2,  then we have
	\begin{equation}
	\norm{[H^{(t)}]^{-1}\nabla F(x^{(t)}) - p^{(t)}}_M  \leq \epsilon_1 \norm{x^{(t)} - x^*}_M.
	\end{equation}
\end{lemma}
\begin{proof}
	\begin{align*}
	\norm{[H^{(t)}]^{-1}\nabla F(x^{(t)}) - p^{(t)}}_M =& \norm{M^{1/2} [H^{(t)}]^{-1}\left(\nabla F(x^{(t)}) -H^{(t)} p^{(t)}\right)}\\
	\overset{\eqref{eq:inexact_p}}{\leq}&\epsilon_1(1+\epsilon_0)^{-1}\kappa^{-3/2}\norm{M^{1/2}} \norm{[H^{(t)}]^{-1}} \norm{\nabla F(x^{(t)})}\\
	\overset{\eqref{eq:prec_cond}}{\leq}&\epsilon_1 \kappa^{-3/2} \norm{M^{1/2}}\norm{[\nabla^2 F(x^{(t)})]^{-1}} \norm{\nabla F(x^{(t)})}\\
	\leq&\epsilon_1\kappa^{-1/2} \norm{M^{1/2}} \norm{x^{(t)} - x^*}\\
	\leq&\epsilon_1\norm{x^{(t)} - x^*}_M,
	\end{align*}
	where the last two inequalities follow from the assumptions that $F(x)$ is $L$-smooth and $\mu$-strongly convex.
\end{proof}

\begin{lemma}
	\label{lem:un_lip}
	There exists a sufficient small value $\gamma$, $\nu = o(1)$, such that when $\|x^{(t)} - x^*\|_M\leq \gamma$, the sequence $\{x^{(t)}\}$ of Algorithm~\ref{alg:app_newton} satisfies 
	\begin{equation*}
	\norm{x^{(t+1)} - x^*}_M \leq \left(\epsilon_0 + \epsilon_1 + 2\nu\mu^{-1}+2\left(2\nu^{1/2}\mu^{-1/2} + \nu\mu^{-1}\right)(\nu\mu^{-1}+1)\right)\norm{x^{(t)} - x^*}_M.
	\end{equation*}
\end{lemma}
\begin{proof}
	Because $\nabla^2 F(x)$ is continuous around $x^*$, then existing a sufficient small value $\gamma$ such that if $\norm{x^{(t)} - x^*}_M \leq \gamma$, then it holds that \citep{ortega1970iterative}
	\begin{equation}
	\label{eq:nu}
	\norm{\nabla^2 F(x^*) - \nabla^2 F(x^{(t)})} \leq \nu,
	\end{equation}
	and 
	\begin{equation}
	\label{eq:high_order}
	\norm{\nabla F(x^{(t)}) - \nabla F(x^*) - \nabla^2 F(x^*)(x^{(t)} - x^*)}_M \leq \nu \norm{x^{(t)} - x^*}_M.
	\end{equation}
	
	By Lemma~\ref{lem:cond}, we have 
	\begin{align*}
	\norm{M^{1/2}[H^{(t)}]^{-1}M^{1/2}} \leq& 1+\epsilon_0  +\mu^{-1/2}\norm{\Delta}^{1/2}(1+\epsilon_0)\left(2+\mu^{-1/2}\norm{\Delta}^{1/2}\right)\\
	\leq&2+2\mu^{-1/2}\norm{\Delta}^{1/2}\left(2+\mu^{-1/2}\norm{\Delta}^{1/2}\right)\\
	\overset{\eqref{eq:nu}}{\leq}&
	2+2\left(2\nu^{1/2}\mu^{-1/2} + \nu\mu^{-1}\right).
	\end{align*}
	
	Combining with Lemma~\ref{lem:decr}, \ref{lem:cond} and \ref{lem:inexact}, we can obtain that
	\begin{align*}
	\norm{x^{(t+1)} - x^*}_M \leq& \left(\epsilon_0 + \epsilon_1 + 4\nu^{1/2}\mu^{-1/2} + 2\nu\mu^{-1}\right)\norm{x^{(t)} - x^*}_M \\
	&+\norm{M^{1/2}[H^{(t)}]^{-1} \left(\nabla F(x^{(t)}) - \nabla^2 F(x^*) (x^{(t)} - x^*) \right) }\\
	\overset{\eqref{eq:high_order}}{\leq}&
	\left(\epsilon_0 + \epsilon_1 + 4\nu^{1/2}\mu^{-1/2} + 2\nu\mu^{-1}\right)\norm{x^{(t)} - x^*}_M 
	\\&+ \nu \norm{M^{-1}}\norm{M^{1/2}[H^{(t)}]^{-1}M^{1/2}}\norm{x^{(t)} - x^*}_M\\
	\leq&
	\left(\epsilon_0 + \epsilon_1 + 2\nu\mu^{-1}+2\left(2\nu^{1/2}\mu^{-1/2} + \nu\mu^{-1}\right)(\nu\mu^{-1}+1)\right)\norm{x^{(t)} - x^*}_M.
	\end{align*}

	From above equation, we can observe that if $\epsilon_0+\epsilon_1 < 1$ and $\nu$ is sufficiently small which can be guaranteed by choosing proper $\gamma$, then we have $\norm{x^{(t+1)} - x^*}_M \leq \norm{x^{(t)} - x^*}_M \leq \gamma$.
\end{proof}

\begin{lemma}
	\label{lem:lip}
	Let the Hessian of $F(x)$ be $\hL$-Lipschitz continuous and the $x^{(t)}$ satisfy $\norm{x^{(t)} - x^*}_M \leq \mu^{3/2} \hL^{-1}$. Then the sequence $\{x^{(t)}\}$ of Algorithm~\ref{alg:app_newton} satisfies 
	\begin{equation*}
	\norm{x^{(t+1)} - x^*}_M  
	\leq 
	\left(\epsilon_0+\epsilon_1\right) \norm{x^{(t)} - x^*}_M  + 7 \mu^{-3/4} \hL^{1/2} \norm{x^{(t)} - x^*}_M^{3/2}.
	\end{equation*}
\end{lemma}
\begin{proof}
	By Taylor's expansion at $x^*$, we have
	\begin{equation*}
	\nabla F(x^{(t)}) - \nabla^2 F(x^*) (x^{(t)} - x^*) = \int_{0}^{1} \nabla^2 F\left(x^* + s(x^{(t)} - x^*)\right) - \nabla^2 F(x^*)\; ds \cdot(x^{(t)} - x^*).
	\end{equation*}
	Thus, we can obtain that
	\begin{align*}
	&\norm{M^{1/2}[H^{(t)}]^{-1} \left(\nabla F(x^{(t)}) - \nabla^2 F(x^*) (x^{(t)} - x^*) \right) }\\
	=&\norm{M^{1/2}[H^{(t)}]^{-1}M^{1/2}  \int_{0}^{1} M^{-1/2} \left(\nabla^2 F\left(x^* + s(x^{(t)} - x^*)\right) - \nabla^2 F(x^*)\right)M^{-1/2} \; ds \cdot M^{1/2}(x^{(t)} - x^*)}\\
	\leq&\underbrace{\norm{M^{1/2}[H^{(t)}]^{-1}M^{1/2}}}_{T_1}\cdot\underbrace{\norm{\int_{1}^{1} \left(M^{-1/2} \left(\nabla^2 F(x^* + s(x^{(t)} - x^*))\right) M^{-1/2} - I\right)\; ds}}_{T_2} \cdot\norm{x^{(t)} - x^*}_M .
	\end{align*}
	Next, we will bound the value of $T_1$ and $T_2$. 
	By Lemma~\ref{lem:cond}, we have
	\begin{align*}
	\norm{M^{1/2}[H^{(t)}]^{-1}M^{1/2}} 
	\leq
	2+2\mu^{-1/2}\norm{\Delta}^{1/2}\left(2+\mu^{-1/2}\norm{\Delta}^{1/2}\right).
	\end{align*}
	with $\Delta = \nabla^2F(x^*) - \nabla^2F(x^{(t)})$. 
	By the assumption that $\nabla^{2}F(x)$ is $\hL$-Lipschitz continuous, then we have
	\begin{align}
	\mu^{-1/2}\norm{\Delta}^{1/2} \left(2+ \mu^{-1/2}\norm{\Delta}^{1/2}\right) \leq & 
	\mu^{-1/2}\hL^{1/2}\norm{x^{(t)} - x^*}^{1/2} \left(2 + \mu^{-1/2}\hL^{1/2} \norm{x^{(t)} - x^*}^{1/2}\right)\notag\\
	\leq& \hL^{1/2} \mu^{-3/4} \norm{x^{(t)} - x^*}_M^{1/2}\left(2 + \mu^{-3/4}\hL^{1/2}\norm{x^{(t)} - x^*}_M^{1/2}\right) \notag\\
	\leq&3 \mu^{-3/4} \hL^{1/2}\norm{x^{(t)} - x^*}_M^{1/2}, \label{eq:high_order_1}\\
	\leq& 3, \notag
	\end{align}
	where the last two inequalities follow from the condition $\norm{x^{(t)} - x^*}_M \leq \mu^{3/2} \hL^{-1}$. 
	Therefore, we can obtain that
	\begin{align}
	T_1 
	\leq 2+2\mu^{-1/2}\norm{\Delta}^{1/2}\left(2+\mu^{-1/2}\norm{\Delta}^{1/2}\right)
	\leq 8. \label{eq:T_1}
	\end{align}

	Let us represent that $$\nabla^2 F(x^* + s(x^{(t)} - x^*)) = M + \Delta',$$
	then we have
	\begin{align*}
	T_2 =& \norm{\int_{0}^{1} \left(M^{-1/2}(M + \Delta') M^{-1/2} - I\right)  \; ds}\\
	=&
	\norm{\int_{0}^{1} \left(M^{-1/2} \Delta'M^{-1/2}\right)\; ds}\\
	\leq&
	\norm{M^{-1}}\int_{0}^{1}\norm{\Delta'} \;ds\\
	\leq&
	\mu^{-1}\hL\int_{0}^{1} \norm{s(x^{(t)} - x^*)}\;ds\\
	\leq&
	\frac{\mu^{-3/2}\hL}{2}\norm{x^{(t)} - x^*}_M.
	\end{align*}
	
	Therefore, we have
	\begin{align*}
	&\norm{M^{1/2}[H^{(t)}]^{-1} \left(\nabla F(x^{(t)}) - \nabla^2 F(x^*) (x^{(t)} - x^*) \right) } \\
	 \leq& T_1 \cdot T_2 \norm{x^{(t)} - x^*}_M\\
 \overset{\eqref{eq:T_1}}{\leq}&8 \cdot \frac{\mu^{-3/2}\hL}{2}  \norm{x^{(t)} - x^*}_M\\
 \leq& 4 \mu^{-3/2}\hL \norm{x^{(t)} - x^*}_M.
	\end{align*}
	Combining with Lemma~\ref{lem:decr}, \ref{lem:cond} and \ref{lem:inexact}, we can obtain that
	\begin{align*}
	\norm{x^{(t+1)} - x^*}_M 
	\leq& 
	\left(\epsilon_0+\epsilon_1 +2\mu^{-1/2}\norm{\Delta}^{1/2}\left(2+\mu^{-1/2}\norm{\Delta}^{1/2}\right)\right)\norm{x^{(t)} - x^*}_M\\
	&+4 \mu^{-3/2}\hL \norm{x^{(t)} - x^*}_M\norm{x^{(t)} - x^*}_M^2\\
	\overset{\eqref{eq:high_order_1}}{\leq}&
	\left(\epsilon_0+\epsilon_1\right) \norm{x^{(t)} - x^*}_M 
	+ 3 \mu^{-3/4} \hL^{1/2} \norm{x^{(t)} - x^*}_M^{3/2} \\
	&+4 \mu^{-3/2}\hL\norm{x^{(t)} - x^*}_M^2\\
	\leq&
	\left(\epsilon_0+\epsilon_1\right) \norm{x^{(t)} - x^*}_M  + 7 \mu^{-3/4} \hL^{1/2} \norm{x^{(t)} - x^*}_M^{3/2},
	\end{align*}
	where the last two inequality follows from the condition $\norm{x^{(t)} - x^*}_M \leq \mu^{3/2} \hL^{-1}$. 
\end{proof}

\section{Proof of Theorem~\ref{thm:glb}}

For a self-concordant function $F(x)$, if two points $x,y$ satisfy $\|x-y\|_x <1$, where $\norm{v}_x = \norm{[\nabla^2 F(x)]^{-1/2} v}$, we have some useful inequalities:
\begin{enumerate}
	\item Hessian bound:
	\begin{equation}
	\label{eq:conc_H}
	(1 - \|x-y\|_x)^2 \nabla^2 F(y) \preceq \nabla^2F(x) \preceq \frac{1}{(1 - \|x-y\|_x)^2} \nabla^2 F(y)
	\end{equation}
	\item Function value bound:
	\begin{equation}
	\label{eq:conc_val}
	\zeta(\|y-x\|_x) \leq F(y) - F(x) - \nabla F(x)^T (y-x) \leq \zeta^*(\|y-x\|_x),
	\end{equation}
	where $\zeta(\alpha) = \alpha - \log(1+\alpha)$ and $\zeta^*(\alpha) = -\alpha - \log(1 - \alpha)$.
\end{enumerate}

This section, we will prove the convergence rate of damped approximate Newton method. 
First, we will show the case that $V(x)$ is smaller than a threshold which is mainly determined by how well the Hessian is approximated. 
In this case, the  step size $s = 1$ will satisfy the exit condition of line search.
Then, we will provide the convergence  analysis when $V(x)$ is larger than the threshold where the step size $s$ should be chosen by the line search.

Before proving the convergence analysis, 
we first define some new notation and clarify their relation.
Let us denote
\begin{align}
V(x^{(t)}) = & \norm{[\nabla^2 F(x^{(t)})]^{-1/2}\nabla F(x^{(t)})}, \label{eq:V}\\
\tV(x^{(t)}) =& \norm{[H^{(t)}]^{-1/2}\nabla F(x^{(t)})}, \label{eq:tV}
\end{align}
and 
\begin{equation}
\hV(x^{(t)}) = \left(\nabla^T F(x^{(t)}) p^{(t)}\right)^{1/2} \label{eq:hV}.
\end{equation}
\begin{lemma}
	\label{lem:rel}
	Let the approximate Hessian satisfy Eqn.~\eqref{eq:prec_cond} and  the descent direction $p^{(t)}$ satisfy Eqn.~\eqref{eq:inexact_p}. Then it holds that 
	\begin{equation*}
	\hV^2(x^{(t)}) \geq \left(1 - \epsilon_1\kappa^{-1}\cdot\left(\frac{1+\epsilon_0}{1-\epsilon_0}\right)^{1/2}\right) \cdot \tV^2(x^{(t)}),
	\end{equation*}
	and
	\begin{equation*}
	\norm{p^{(t)}}^2_{x^{(t)}} 
	\leq 
	(1+\epsilon_0)
	\left(1 + \epsilon_1\kappa^{-1}\cdot\left(\frac{1+\epsilon_0}{1-\epsilon_0}\right)^{1/2}\right)^2 \cdot \tV^2(x^{(t)}).
	\end{equation*}
\end{lemma}
\begin{proof}
	First, we have
	\begin{align*}
	\nabla^T F(x^{(t)}) p^{(t)} 
	=& 
	\nabla^T F(x^{(t)}) [H^{(t)}]^{-1}  \nabla F(x^{(t)}) 
	+ 
	\nabla^T F(x^{(t)})[H^{(t)}]^{-1}\left( [H^{(t)}]p^{(t)} -   \nabla F(x^{(t)})\right)
	\\
	\overset{\eqref{eq:inexact_p}}{\geq}&
	\nabla^T F(x^{(t)}) [H^{(t)}]^{-1}  \nabla F(x^{(t)}) \\
	&-
	\left(\nabla^T F(x^{(t)}) [H^{(t)}]^{-1}  \nabla F(x^{(t)}) \right)^{1/2} \norm{H^{-1/2}} \kappa^{-3/2}\epsilon_1\norm{\nabla F(x^{(t)})}
	\\
	\geq&
	\nabla^T F(x^{(t)}) [H^{(t)}]^{-1}  \nabla F(x^{(t)}) \\
	&-
	\kappa^{-3/2}\epsilon_1 \nabla^T F(x^{(t)}) [H^{(t)}]^{-1}  \nabla F(x^{(t)}) \norm{H^{-1/2}}\norm{H^{1/2}}\\
	\overset{\eqref{eq:prec_cond}}{\geq}&
	\left(1 - \epsilon_1\kappa^{-1}\cdot\left(\frac{1+\epsilon_0}{1-\epsilon_0}\right)^{1/2}\right) \cdot \tV^2(x^{(t)}).
	\end{align*}
	Similarly, we can obtain that
	\begin{equation}
	\label{eq:nab_nab}
	\nabla^T F(x^{(t)}) p^{(t)} 
	\leq 
	\left(1 + \epsilon_1\kappa^{-1}\cdot\left(\frac{1+\epsilon_0}{1-\epsilon_0}\right)^{1/2}\right) 
	\cdot 
	\tV^2(x^{(t)}).
	\end{equation}
	
	By the condition~\eqref{eq:prec_cond}, we can obtain that
	\begin{equation}
	\label{eq:up}
	\norm{p^{(t)}}_{x^{(t)}}^2 \leq (1+\epsilon_0) [p^{(t)}]^T [H^{(t)}] p^{(t)}.
	\end{equation}
	Furthermore, we have
	\begin{align*}
	[p^{(t)}]^T [H^{(t)}] p^{(t)} 
	=& 
	[p^{(t)}]^T \left(\nabla F(x^{(t)}) + [H^{(t)}] p^{(t)} - \nabla F(x^{(t)})\right)\\
	\leq&[p^{(t)}]^T \nabla F(x^{(t)}) 
	+ 
	\norm{p^{(t)}} \norm{[H^{(t)}] p^{(t)} - \nabla F(x^{(t)})}\\
	\overset{\eqref{eq:inexact_p}}{\leq}&
	[p^{(t)}]^T \nabla F(x^{(t)})
	+
	\epsilon_1\kappa^{-3/2}\norm{p^{(t)}}\norm{\nabla F(x^{(t)})}.
	\end{align*}
	Furthermore, we have
	\begin{align*}
	\norm{p^{(t)}}
	\leq&
	\norm{p^{(t)} - [H^{(t)}]^{-1}\nabla F(x^{(t)})} 
	+ 
	\norm{[H^{(t)}]^{-1}\nabla F(x^{(t)})}\\
	\overset{\eqref{eq:inexact_p}}{\leq}& 
	\epsilon_1\kappa^{-3/2}\norm{[H^{(t)}]^{-1}} \norm{\nabla F(x^{(t)})}
	+ 
	\norm{[H^{(t)}]^{-1/2}}\norm{[H^{(t)}]^{-1/2}\nabla F(x^{(t)})}
	\\
	\leq&
	\left(
	\epsilon_1\kappa^{-3/2}\norm{[H^{(t)}]^{-1}}\norm{[H^{(t)}]^{1/2}} + \norm{[H^{(t)}]^{-1/2}}
	\right)
	\norm{[H^{(t)}]^{-1/2}\nabla F(x^{(t)})}
	\\
	\leq&
	\left(
	1 + \epsilon_1\kappa^{-1}\cdot\left(\frac{1+\epsilon_0}{1-\epsilon_0}\right)^{1/2}
	\right) 
	\norm{[H^{(t)}]^{-1/2}}
	\norm{[H^{(t)}]^{-1/2}\nabla F(x^{(t)})}
	\end{align*}
	and
	\begin{align*}
	\norm{\nabla F(x^{(t)})} 
	\leq
	\norm{[H^{(t)}]^{1/2}}
	\norm{[H^{(t)}]^{-1/2}\nabla F(x^{(t)})}.
	\end{align*}
	Thus, we can obtain that
	\begin{align*}
	\norm{p^{(t)}} \norm{\nabla F(x^{(t)})} 
	\leq& 
	\left(
	1 + \epsilon_1\kappa^{-1}\cdot\left(\frac{1+\epsilon_0}{1-\epsilon_0}\right)^{1/2}
	\right)
	\norm{[H^{(t)}]^{-1/2}}
	\norm{[H^{(t)}]^{1/2}}
	\norm{[H^{(t)}]^{-1/2}\nabla F(x^{(t)})}^2
	\\
	\leq&
	\kappa^{1/2}\left(\frac{1+\epsilon_0}{1-\epsilon_0}\right)^{1/2} 
	\left(
	1 + \epsilon_1\kappa^{-1}\cdot\left(\frac{1+\epsilon_0}{1-\epsilon_0}\right)^{1/2}
	\right)
	\tV^2(x^{(t)}).
	\end{align*}
	Therefore, we can obtain that
	\begin{align*}
	[p^{(t)}]^T [H^{(t)}] p^{(t)}  
	\leq&
	[p^{(t)}]^T \nabla F(x^{(t)})
	+
	\epsilon_1\kappa^{-3/2}\norm{p^{(t)}}\norm{\nabla F(x^{(t)})}
	\\
	\leq&
	[p^{(t)}]^T \nabla F(x^{(t)})
	+
	\epsilon_1\kappa^{-1}
	\left(\frac{1+\epsilon_0}{1-\epsilon_0}\right)^{1/2} 
	\left(
	1 + \epsilon_1\kappa^{-1}\cdot\left(\frac{1+\epsilon_0}{1-\epsilon_0}\right)^{1/2}
	\right)
	\tV^2(x^{(t)})\\
	\overset{\eqref{eq:nab_nab}}{\leq}&
	\left(1 + \epsilon_1\kappa^{-1}\cdot\left(\frac{1+\epsilon_0}{1-\epsilon_0}\right)^{1/2}\right) ^2
	\tV^2(x^{(t)}).
	\end{align*}
	Combining Eqn.~\eqref{eq:up}, we can obtain
	\begin{equation*}
	\norm{p^{(t)}}^2_{x^{(t)}} 
	\leq 
	(1+\epsilon_0)
	\left(1 + \epsilon_1\kappa^{-1}\cdot\left(\frac{1+\epsilon_0}{1-\epsilon_0}\right)^{1/2}\right)^2 \cdot \tV^2(x^{(t)}).
	\end{equation*}
\end{proof}

Now, we begin to prove the case that $	V(x^{(t)}) \leq \frac{1-\epsilon_0 - 2\epsilon_1\kappa^{-1}}{12}$ and the step size $s = 1$ is sufficient.
\begin{lemma}
	\label{lem:local}
	Let the descent direction $p^{(t)}$ satisfy Eqn.~\eqref{eq:inexact_p} and $V(x^{(t)})$ satisfy 
	\[
	V(x^{(t)}) \leq \frac{1-\epsilon_0 - 2\epsilon_1\kappa^{-1}}{12}.
	\]
	Then the approximate Newton with backtrack line search (Algorithm~\ref{alg:damp_app_newton}) has the following convergence property
	\begin{equation*}
	V(x^{(t+1)}) \leq \frac{1+\epsilon_0+2\epsilon_1\kappa^{-1}}{2} V(x^{(t)}).
	\end{equation*}
\end{lemma}
\begin{proof}
	Then we have
	\begin{align*}
	V(x^{(t+1)})=& \norm{[\nabla^2F(x^{(t+1)})]^{-1/2}\nabla F(x^{(t+1)})}\\
	\overset{\eqref{eq:conc_H}}{\leq}& 
	\frac{1}{1 - \norm{p^{(t)}}_x}\norm{[\nabla^2F(x^{(t)})]^{-1/2}\nabla F(x^{(t+1)})}
	\end{align*}
	By Taylor's expansion of $\nabla F(x^{(t+1)})$ at point $x^{(t)}$, we have
	\begin{align*}
	&\norm{[\nabla^2F(x^{(t)})]^{-1/2}\nabla F(x^{(t+1)})} \\
	= &\norm{[\nabla^2F(x^{(t)})]^{-1/2}\left(\nabla F(x^{(t)})+\nabla^2 F(x^{(t)})(-p^{(t)}) + \int_0^1 [\nabla^2 F(x^{(t)} + sp^{(t)}) - \nabla^2 F(x^{(t)})](-p^{(t)})ds\right)}\\
	\leq&\underbrace{\norm{\left(I - [\nabla^2F(x^{(t)})]^{1/2}[H^{(t)}]^{-1}[\nabla^2F(x^{(t)})]^{1/2}\right)[\nabla^2F(x^{(t)})]^{-1/2}\nabla F(x^{(t)})}}_{T_1}\\
	&+\underbrace{\norm{[\nabla^2F(x^{(t)})]^{1/2}[H^{(t)}]^{-1}[\nabla^2F(x^{(t)})]^{1/2}}
		\cdot 
		\norm{[\nabla^2F(x^{(t)})]^{-1/2}}
		\cdot \norm{\nabla F(x^{(t)}) - H^{(t)} p^{(t)}}}_{T_2} \\
	&+\underbrace{\norm{\int_0^1 \left([\nabla^2F(x^{(t)})]^{-1/2}\nabla^2 F(x^{(t)} - sp^{(t)}) [\nabla^2F(x^{(t)})]^{-1/2} - I\right)ds \cdot [\nabla^2F(x^{(t)})]^{1/2} p^{(t)}}}_{T_3}
	\end{align*}
	We are going to bound the above terms. 
	First, by the assumption~\eqref{eq:prec_cond}, we have 
	\begin{equation*}
	\norm{I - [\nabla^2F(x^{(t)})]^{1/2}[H^{(t)}]^{-1}[\nabla^2F(x^{(t)})]^{1/2}} \leq \epsilon_0.
	\end{equation*}
	Combining the definition of $V(x)$, we can obtain
	\begin{align*}
	T_1 \leq& \norm{I - [\nabla^2F(x^{(t)})]^{1/2}[H^{(t)}]^{-1}[\nabla^2F(x^{(t)})]^{1/2}}\cdot\norm{[\nabla^2F(x^{(t)})]^{-1/2}\nabla F(x^{(t)})} \\
	\leq& \epsilon_0 V(x^{(t)}). 
	\end{align*}
	
	Also by the condition~\eqref{eq:prec_cond}, we have
	\begin{equation*}
	\norm{[\nabla^2F(x^{(t)})]^{1/2}[H^{(t)}]^{-1}[\nabla^2F(x^{(t)})]^{1/2}} \leq (1+\epsilon_0).
	\end{equation*}
	Combining the condition \eqref{eq:inexact_p} and the definition of $V^{(t)}$, we can obtain that
	\begin{equation*}
	T_2 \leq (1+\epsilon_0) \mu^{-1/2} \frac{\epsilon_1}{\kappa^{3/2}}\norm{\nabla F(x^{(t)})} \leq \frac{(1+\epsilon_0)\epsilon_1}{\kappa} V(x^{(t)}) \leq \frac{2\epsilon_1}{\kappa} V(x^{(t)}).
	\end{equation*}
	We also have
	\begin{align*}
	T_3 \leq&\norm{\int_0^1 \left([\nabla^2F(x^{(t)})]^{-1/2}\nabla^2 F(x^{(t)} - sp^{(t)}) [\nabla^2F(x^{(t)})]^{-1/2} - I\right)ds }\cdot\norm{p^{(t)}}_x\\
	\overset{\eqref{eq:conc_H}}{\leq}&
	\left|\int_{0}^{1}\left(\frac{1}{(1 - s\norm{p^{(t)}}_x)^2} - 1\right) ds\right|
	\cdot 
	\norm{p^{(t)}}_x\\
	=&\frac{\norm{p^{(t)}}_x}{1-\norm{p^{(t)}}_x} 
	\cdot 
	\norm{p^{(t)}}_x.
	\end{align*}
	Next, we will bound the value of $\norm{p^{(t)}}_x$. 
	We have
	\begin{align*}
	\norm{p^{(t)}}_x =& \norm{[\nabla^2 F(x^{(t)})]^{1/2} p^{(t)}}\\
	=&\norm{[\nabla^2 F(x^{(t)})]^{1/2} [H^{(t)}]^{-1} \nabla F(x^{(t)}) - [\nabla^2 F(x^{(t)})]^{1/2} [H^{(t)}]^{-1} \left( \nabla F(x^{(t)}) - H^{(t)} p^{(t)}\right)}\\
	\leq&
	\norm{[\nabla^2 F(x^{(t)})]^{1/2} [H^{(t)}]^{-1} [\nabla^2 F(x^{(t)})]^{1/2} \cdot [\nabla^2 F(x^{(t)})]^{-1/2} \nabla F(x^{(t)})}\\
	&+
	\norm{[\nabla^2F(x^{(t)})]^{1/2}[H^{(t)}]^{-1}[\nabla^2F(x^{(t)})]^{1/2}}
	\cdot 
	\norm{[\nabla^2F(x^{(t)})]^{-1/2}}
	\cdot \norm{\nabla F(x^{(t)}) - H^{(t)} p^{(t)}}\\
	=& 
	\norm{[\nabla^2 F(x^{(t)})]^{1/2} [H^{(t)}]^{-1} [\nabla^2 F(x^{(t)})]^{1/2} \cdot [\nabla^2 F(x^{(t)})]^{-1/2} \nabla F(x^{(t)})} + T_2\\
	\leq&
	(1+\epsilon_0) V(x^{(t)}) + \frac{2\epsilon_1}{\kappa} V(x^{(t)}).
	\end{align*}
	Combining above results, we can obtain that
	\begin{align*}
	V(x^{(t+1)}) 
	\leq& 
	\frac{1}{1 - \norm{p^{(t)}}_x} (T_1 + T_2+T_3)\\
	\leq&
	\frac{(\epsilon_0+2\epsilon_1\kappa^{-1}) V(x^{(t)})}{1 - (1+\epsilon_0+2\epsilon_1\kappa^{-1}) V(x^{(t)})} + \frac{(1+\epsilon_0+2\epsilon_1\kappa^{-1})^2 V^2(x^{(t)})}{(1-(1+\epsilon_0+2\epsilon_1\kappa^{-1})V(x^{(t)}))^2}
	\end{align*}
	
	If $V(x^{(t)})$ satisfies that 
	\begin{align}
	V(x^{(t)}) \leq&
	\frac{1 - (\epsilon_0+2\epsilon_1\kappa^{-1})^2}{(1+\epsilon_0+2\epsilon_1\kappa^{-1})^2\left(2+\epsilon_0 +2\epsilon_1\kappa^{-1}+ \sqrt{(2+\epsilon_0+2\epsilon_1\kappa^{-1})^2 - 1 +(\epsilon_0+ 2\epsilon_1\kappa^{-1})^2}\right)} \notag \\
	\leq& 
	\frac{1-\epsilon_0 - 2\epsilon_1\kappa^{-1}}{12}, \label{eq:vx}
	\end{align}
	we have
	\begin{equation*}
	V(x^{(t+1)}) \leq \frac{1+\epsilon_0+2\epsilon_1\kappa^{-1}}{2} V(x^{(t)}).
	\end{equation*}
\end{proof}

Now we begin to analyze the phase that line search should be applied to find a step size $s<1$. 
This phase is commonly commonly referred as \emph{damped phase}.    
\begin{lemma}
	\label{lem:damp}
	Let the approximate Hessian satisfy Eqn.~\eqref{eq:prec_cond} and  the descent direction $p^{(t)}$ satisfy Eqn.~\eqref{eq:inexact_p}. 
	If it holds that 
	\[
	\tV(x) \geq \frac{\sqrt{(1-\epsilon_0)} (1-\epsilon_0 - 2\epsilon_1\kappa^{-1}) }{12},
	\]
	then Algorithm~\ref{alg:damp_app_newton} has the following convergence property
	\begin{equation*}
	F(x^{(t+1)}) \leq 
	F(x^{(t)}) -\alpha \beta\cdot \frac{\rho^2 \tV^2(x^{(t)})}{1+\rho \tV(x^{(t)})},
	\end{equation*}
	where $\rho$ is defined as
	\begin{equation*}
	\rho =\frac{(1-\varphi)^{1/2} }{(1+\epsilon_0)^{1/2} (1+\varphi)}, \quad\mbox{with}\quad \varphi = \epsilon_1\kappa^{-1}\cdot\left(\frac{1+\epsilon_0}{1-\epsilon_0}\right)^{1/2}.
	\end{equation*}
\end{lemma}
\begin{proof}
	By the update rule, we can obtain that
	\begin{align*}
	F(x^{(t+1)}) \overset{\eqref{eq:conc_val}}{\leq}& F(x^{(t)}) - s \nabla F(x^{(t)})^Tp^{(t)} + \zeta^*\left(s \norm{p^{(t)}}_{x^{(t)}}\right) \\
	=&F(x^{(t)}) - s \hV^2(x^{(t)}) - s \norm{p^{(t)}}_{x^{(t)}} - \log\left(1 - s \norm{p^{(t)}}_{x^{(t)}})\right),
	\end{align*}
	with $0\leq s < 1/\tV(x^{(t)}) $.

	Letting us define $\hat{s}$ as
	\begin{equation*}
	\hat{s} = \frac{\hV^2(x^{(t)})}{\left(\hV^2(x^{(t)}) + \norm{p^{(t)}}_{x^{(t)}}\right) \norm{p^{(t)}}_{x^{(t)}}}.
	\end{equation*}
	We can use this bound to show the backtracking line search always results in a step size $s \geq \beta \hat{s}$. 
	Furthermore, we can obtain that 
	\begin{align*}
	F(x^{(t+1)}) \leq& 
	F(x^{(t)}) - \frac{\hV^2(x^{(t)})}{\norm{p^{(t)}}_{x^{(t)}}} - \log\left(\frac{\norm{p^{(t)}}_{x^{(t)}}}{\hV^2(x^{(t)}) + \norm{p^{(t)}}_{x^{(t)}}}\right)\\
	=&  F(x^{(t)}) -  \frac{\hV^2(x^{(t)})}{\norm{p^{(t)}}_{x^{(t)}}} + \log\left(1 +\frac{\hV^2(x^{(t)})}{\norm{p^{(t)}}_{x^{(t)}}}\right)\\
	\leq&F(x^{(t)}) - \frac{\left(\hV^2(x^{(t)})/\norm{p^{(t)}}_{x^{(t)}})\right)^2}{2\left(1+\hV^2(x^{(t)})/\norm{p^{(t)}}_{x^{(t)}}\right)},\\
	=&
	F(x^{(t)}) - \frac{1}{2}\cdot \hat{s} \hV^2(x^{(t)})\\
	\leq& F(x^{(t)}) -\alpha \cdot \hat{s} \hV^2(x^{(t)})
	\end{align*}
	where the second inequality follows form the fact that it holds  for $a>0$ that
	\begin{equation*}
	-a + \log(1+a) + \frac{a^2}{2(1+a)} \leq 0.
	\end{equation*}
	The last inequality is because $\alpha < 1/2$.
	Since we obtain that $F(x^{(t+1)}) \leq  F(x^{(t)}) -\alpha \cdot \hat{s} \hV^2(x^{(t)})$, we show the exit condition of the line search has satisfied.
	Furthermore, the exit condition holds when  the step size satisfies $s \geq \beta \hat{s}$.
	Thus, we can obtain that
	\begin{align*}
	F(x^{(t+1)}) \leq& F(x^{(t)}) -\alpha \beta\cdot \hat{s} \hV^2(x^{(t)}).
	\end{align*}
	Next, we will bound the value of $\hat{s} \hV^2(x^{(t)})$.
	By the definition of $\hat{s}$, we can obtain that 
	\begin{equation*}
	\hat{s} \hV^2(x^{(t)}) 
	= 
	\frac{\left(\hV^2(x^{(t)})/\norm{p^{(t)}}_{x^{(t)}}\right)^2}{\left(1+\hV^2(x^{(t)})/\norm{p^{(t)}}_{x^{(t)}}\right)}. 
	\end{equation*} 
	By Lemma~\ref{lem:rel}, we have
	\begin{equation*}
	\frac{\hV(x^{(t)})}{\norm{p^{(t)}}_{x^{(t)}}} 
	\geq 
	\frac{(1-\varphi)^{1/2} \tV(x^{(t)})}{(1+\epsilon_0)^{1/2} (1+\varphi)\tV(x^{(t)})}
	\\
	=\frac{(1-\varphi)^{1/2} }{(1+\epsilon_0)^{1/2} (1+\varphi)},
	\end{equation*}
	where $\varphi = \epsilon_1\kappa^{-1}\cdot\left(\frac{1+\epsilon_0}{1-\epsilon_0}\right)^{1/2}$. 
	Furthermore, we have 
	\begin{align*}
	\hat{s} \hV^2(x^{(t)}) 
	=& \frac{\left(\hV^2(x^{(t)})/\norm{p^{(t)}}_{x^{(t)}}\right)^2}{\left(1+\hV^2(x^{(t)})/\norm{p^{(t)}}_{x^{(t)}}\right)}\\
	\geq&
	\frac{\left(\frac{(1-\varphi)^{1/2} }{(1+\epsilon_0)^{1/2} (1+\varphi)}\hV(x^{(t)})\right)^2}{1+ \frac{(1-\varphi)^{1/2} }{(1+\epsilon_0)^{1/2} (1+\varphi)}\hV(x^{(t)})}
	\\
	\geq&
	\frac{\left(\frac{(1-\varphi) }{(1+\epsilon_0)^{1/2} (1+\varphi)}\tV(x^{(t)})\right)^2}{1+ \frac{(1-\varphi) }{(1+\epsilon_0)^{1/2} (1+\varphi)}\tV(x^{(t)})}
	\end{align*}
	where the last inequality follows from Lemma~\ref{lem:rel}.
	
	Letting us denote $\rho = \frac{(1-\varphi)^{1/2} }{(1+\epsilon_0)^{1/2} (1+\varphi)}$, then we have
	\begin{align*}
	F(x^{(t+1)}) 
	\leq& 
	F(x^{(t)}) -\alpha \beta\cdot \hat{s} \hV^2(x^{(t)}) 
	\\
	\leq& 
	F(x^{(t)}) -\alpha \beta\cdot \frac{\rho^2 \tV^2(x^{(t)})}{1+\rho \tV(x^{(t)})}
	\end{align*}
	
	By the Condition~\eqref{eq:prec_cond}, we have
	\begin{equation}
	\label{eq:vx_1}
	\frac{1}{1-\epsilon_0} \tV^2(x^{(t)}) \geq V^2(x^{(t)}).
	\end{equation}
	Thus, we can obtain that if $\tV(x) \leq \frac{(1-\epsilon_0)^{1/2} 1-\epsilon_0 - 2\epsilon_1\kappa^{-1} }{12}$, then it holds that $V(x) \leq \frac{1-\epsilon_0}{12}$.
	Therefore, we can obtain that when $\tV(x) \geq \frac{(1-\epsilon_0)^{1/2} 1-\epsilon_0 - 2\epsilon_1\kappa^{-1} }{12}$, it holds that
	\begin{equation*}
	F(x^{(t+1)}) \leq 
	F(x^{(t)}) -\alpha \beta\cdot \frac{\rho^2 \tV^2(x^{(t)})}{1+\rho \tV(x^{(t)})}.
	\end{equation*}
\end{proof}

Combining Lemma~\ref{lem:local} and \ref{lem:damp}, we can obtain the global convergence rate of approximate Newton with backtracking line search.
\begin{proof}{\bf of Theorem~\ref{thm:glb}}
	Let us denote 
	\begin{equation*}
	\eta = \alpha \beta\cdot \frac{\rho^2 \left(\frac{\sqrt{(1-\epsilon_0)} (1-\epsilon_0 - 2\epsilon_1\kappa^{-1}) }{12}\right)^2}{1+\rho \frac{\sqrt{(1-\epsilon_0)} (1-\epsilon_0 - 2\epsilon_1\kappa^{-1}) }{12}} = \alpha\beta \frac{(1-\epsilon_0)\rho^2(1-\epsilon_0 - 2\epsilon_1\kappa^{-1})^2}{144 + 12 \rho \sqrt{(1-\epsilon_0)} (1-\epsilon_0 - 2\epsilon_1\kappa^{-1})}.
	\end{equation*}
	By Lemma~\ref{lem:damp}, we can obtain that it takes at most 
	\begin{equation*}
	\frac{F(x^{(0)}) - F(x^*)}{\eta}
	\end{equation*}
	steps in the damped phase
	because of $F(x^{(t+1)}) - F(x^{(t)}) \leq -\eta$ when $\tV(x) \geq \frac{\sqrt{(1-\epsilon_0)} (1-\epsilon_0 - 2\epsilon_1\kappa^{-1}) }{12}$.
	
	If it holds that  $\tV(x) \leq \frac{\sqrt{(1-\epsilon_0)} (1-\epsilon_0 - 2\epsilon_1\kappa^{-1}) }{12}$, then we have $V(x^{(t)}) \leq \frac{1-\epsilon_0 - 2\epsilon_1\kappa^{-1}}{12}$. 
	By Lemma~\ref{lem:local}, we have
	\begin{equation*}
	V(x^{(t+k)}) \leq \left(\frac{1+\epsilon_0 + 2\epsilon_1\kappa^{-1}}{2}\right)^k \frac{1-\epsilon_0-2\epsilon_1\kappa^{-1}}{12}
	\end{equation*}
	Furthermore, the self-concordance of $F(x)$ implies that
	\begin{equation*}
	F(x^{(t+k)}) - F(x^*) \leq V(x^{(t+k)}) \leq  \left(\frac{1+\epsilon_0 + 2\epsilon_1\kappa^{-1}}{2}\right)^k \frac{1-\epsilon_0 - 2\epsilon_1\kappa^{-1}}{12}.
	\end{equation*} 
	To make the right hand of above equation less than $\epsilon$, then it will take no more than
	\begin{equation*}
	k = \frac{2}{1-\epsilon_0- 2\epsilon_1\kappa^{-1} }\log\left(\frac{1-\epsilon_0 - 2\epsilon_1\kappa^{-1}}{12\epsilon}\right)
	\end{equation*}
	iterations.
	
	Therefore,  the total complexity of approximate Newton method with backtracking line search to achieve an $\epsilon$-suboptimality   is at most 
	\begin{equation*}
	\frac{F(x^{(0)}) - F(x^*)}{\eta} + \frac{2}{1-\epsilon_0- 2\epsilon_1\kappa^{-1} }\log\left(\frac{1-\epsilon_0 - 2\epsilon_1\kappa^{-1}}{12\epsilon}\right).
	\end{equation*}
\end{proof}

\section{Proofs of Section~\ref{sec:ske_newton}}

\begin{proof}{\bf{of Theorem~\ref{thm:Sketch_newton}}}
If $S$ is an $\epsilon_0$-subspace embedding matrix w.r.t.\ $B(x^{(t)})$, then we have
\begin{equation}
(1-\epsilon_0) \nabla^{2}F(x^{(t)})\preceq [B(x^{(t)})]^{T}S^{T}SB(x^{(t)}) \preceq (1+\epsilon_0) \nabla^{2}F(x^{(t)}).\label{eq: preceq_1}
\end{equation}
By simple transformation and omitting $\epsilon_0^{2}$, Eqn.~\eqref{eq: preceq_1} can be transformed into
\[
(1-\epsilon_0) [B(x^{(t)})]^{T}S^{T}S\nabla^{2}B(x^{(t)}) \preceq \nabla^{2}F(x^{(t)}) \preceq (1+\epsilon_0) [B(x^{(t)})]^{T}S^{T}SB(x^{(t)}).
\]
The convergence rate can be derived directly from Theorem~\ref{thm:univ_frm} and \ref{thm:glb}.
\end{proof}	

\begin{proof}{\bf{of Corollary~\ref{cor:suplin}}}
If $\nabla^2F(x)$ is not Lipschitz continuous, then we have
\begin{align*}
\limsup_{t \to \infty} \frac{\|x^{(t+1)} - x^*\|_{M}}{\|x^{(t)} - x^* \|_{M}} &= \limsup_{t \to \infty} \left(\epsilon_0(t)+ \nu(t)\kappa\mu^{-1}(2\mu^{1/2} + 2\kappa^{-1/2} +\nu(t))\right) \\
& = \limsup_{t \to \infty} \left(\frac{1}{\log(1+t)}+ \nu(t)\kappa\mu^{-1}(2\mu^{1/2} + 2\kappa^{-1/2} +\nu(t))\right) \\ 
&= 0,
\end{align*}
where $\nu(t) \to 0$ is because $\|\nabla^2F(x^{(t)}) - \nabla^2F(x^*)\|\to 0$ as $x^{(t)}$ approaches $x^*$.

If $\nabla^2F(x)$ is Lipschitz continuous, then we have
\begin{align*}
\limsup_{t \to \infty} \frac{\|x^{(t+1)} - x^*\|_{M}}{\|x^{(t)} - x^* \|_{M}} 
\leq& 
\limsup_{t \to \infty} \left(\epsilon_0(t) + 7 \mu^{-3/4} \hL^{1/2} \norm{x^{(t)} - x^*}_M^{1/2} \right) \\
=& 
\limsup_{t \to \infty} \left(\frac{1}{\log(1+t)} + 7 \mu^{-3/4} \hL^{1/2} \norm{x^{(t)} - x^*}_M^{1/2}\right) \\
=& 0.
\end{align*}
\end{proof}
\section{Proofs of theorems of Section~\ref{sec:sub_newton}}

\begin{proof}{\bf{of Theorem~\ref{thm:H_subsamp}}}
	Let us denote that
	\begin{equation*}
	X_i = [\nabla^{2}F(x^{(t)})]^{-1/2} \nabla^2f_i(x) [\nabla^{2}F(x^{(t)})]^{-1/2}, \quad\mbox{and}\quad
	Y = \sum_{i\in \SM} X_i
	\end{equation*}
	Because $\nabla^2f_i(x)$ is chosen uniformly, then we have $\EB[Y] = \sum_{i\in \SM} \EB[X_i] = \SM I$. 
	Furthermore, by the Condition~\eqref{eq:k} and \eqref{eq:sigma}, we can obtain that
	\begin{equation*}
	\norm{X_i} \leq \frac{K}{\mu} \quad\mbox{and}\quad \lambda_{\max}(\EB[y]) = \lambda_{\min}(\EB[y]) = |\SM|.
	\end{equation*}
	By Lemma~\ref{lem:matrix_bnd}, we have
	\begin{align*}
	\PB\left(\lambda_{\min}(Y) \leq (1-\epsilon_0) |\SM|\right) \leq d\exp\left(-\frac{\epsilon_0^2|\SM|}{2K/\mu}\right).
	\end{align*}
	Letting us choose $|\SM| = \frac{2K/\mu\log(d/\delta)}{\epsilon_0^2}$, then it holds with probability at least $1-\delta$ that 
	\begin{equation*}
	\lambda_{\min}(Y) \geq 1-\epsilon_0
	\end{equation*}
	which implies that
	\begin{align*}
	&\min_{x\in\RB^{d}} \frac{x^T[\nabla^{2}F(x^{(t)})]^{-1/2} \left(\sum_{i\in \SM}\nabla^2f_i(x)\right) [\nabla^{2}F(x^{(t)})]^{-1/2} x}{\norm{x}^2}\geq (1-\epsilon_0) |\SM|\\
	\Rightarrow& \frac{1}{|\SM|}\sum_{i\in \SM}\nabla^2f_i(x) \succeq (1-\epsilon_0) \nabla^{2}F(x^{(t)}).
	\end{align*}
	By simple transformation and omitting $\epsilon_0^{2}$, the above equation can be represented as
	\begin{equation}
	\label{eq:sub_one}
	\nabla^{2}F(x^{(t)}) \preceq (1+\epsilon_0) H^{(t)}.
	\end{equation}

	Also by Lemma~\ref{lem:matrix_bnd}, we have
	\begin{equation*}
	\PB\left(\lambda_{\max}(Y) \geq (1+\epsilon_0) |\SM|\right) \leq d\exp\left(-\frac{\epsilon_0^2|\SM|}{3K/\mu}\right).
	\end{equation*}
	By the similar proof of above, we can obtain that if we choose $|\SM| = \frac{3K/\mu\log(d/\delta)}{\epsilon_0^2}$, it holds with probability at least $1-\delta$ that 
	\begin{equation*}
	\left(1-\epsilon_0\right) H^{(t)} \preceq \nabla^{2}F(x^{(t)})
	\end{equation*}
	Combining with Eqn.~\eqref{eq:sub_one} and by the union bound of probability, we can obtain that if we choose $|\SM| = \frac{3K/\mu\log(2d/\delta)}{\epsilon_0^2}$, it holds that 
	\[
	(1-\epsilon_0) H^{(t)} \preceq \nabla^{2}F(x^{(t)}) \preceq (1+\epsilon_0) H^{(t)},
	\]
	with probability at least $1-\delta$.
	
	Finally, the local convergence properties of Algorithm~\ref{alg:H_subsamp} can be obtained by Theorem~\ref{thm:univ_frm} and Theorem~\ref{thm:glb}.
\end{proof}

\begin{proof} {\bf{of Theorem~\ref{thm:Reg_subnewton}}}\\
Let us denote that 
\begin{equation*}
X_i = [\nabla^{2}F(x^{(t)})+\xi I]^{-1/2} \left(\nabla^2f_i(x)+\xi I\right) [\nabla^{2}F(x^{(t)}) +\xi I ]^{-1/2}, \quad\mbox{and}\quad
Y = \sum_{i\in \SM} X_i
\end{equation*}
Then we can obtain that
\begin{equation*}
\norm{X_i} \leq \frac{K + \xi}{\mu+\xi}
\end{equation*}
Because $\nabla^2f_i(x)$ is chosen uniformly, then we have $\EB[Y] = \sum_{i\in \SM} \EB[X_i] = \SM I $. 
Hence, we can obtain that
\begin{equation}
\lambda_{\max}(Y) = \lambda_{\min}(Y) = \SM.
\end{equation}

By Lemma~\ref{lem:matrix_bnd}, we have 
\begin{equation*}
\PB\left(\lambda_{\min}(Y) \leq \frac{2}{3} |\SM|\right) \leq d\exp\left(-\frac{|\SM|}{18(K+\xi)/(\mu+\xi)}\right).
\end{equation*}
Letting us choose $|\SM| = \frac{18K\log(d/\delta)}{\xi}$, then it holds with probability at least $1-\delta$ that 
\begin{equation}
\frac{1}{|\SM|}\sum_{i\in \SM}\nabla^2f_i(x) + \xi I\succeq \frac{2}{3} \left(\nabla^{2}F(x^{(t)})+\xi I\right)\succeq \frac{2}{3} \left(1+ \frac{\xi}{L}\right) \nabla^{2}F(x^{(t)}), \label{eq:up_1}
\end{equation}
which implies that 
\begin{equation*}
\nabla^{2}F(x^{(t)}) \preceq \left(1+ \frac{L - 2\xi}{2(L+\xi)}\right) H^{(t)}.
\end{equation*}

Also by Lemma~\ref{lem:matrix_bnd}, we have
\begin{equation*}
\PB\left(\lambda_{\max}(Y) \geq \frac{3}{2} |\SM|\right) \leq d\exp\left(-\frac{|\SM|}{12(K+\xi)/(\mu+\xi)}\right).
\end{equation*}
By the similar proof of above, we can obtain that if we choose $|\SM| = \frac{12K\log(d/\delta)}{\xi}$, it holds with probability at least $1-\delta$ that 
\begin{equation}
\label{eq:down_1}
 \frac{1}{|\SM|}\sum_{i\in \SM}\nabla^2f_i(x) + \xi I \preceq \frac{3}{2}\left(\nabla^{2}F(x^{(t)})+\xi I\right)\preceq \frac{3}{2}\left(1+\frac{\xi}{\mu}\right) \nabla^{2}F(x^{(t)}),
\end{equation}
which implies that 
\begin{equation*}
\left(1 - \frac{3\xi + \mu}{3\alpha + 3\mu}\right) H^{(t)} \preceq\nabla^{2}F(x^{(t)}).
\end{equation*}

Therefore, by choosing $|\SM| = \frac{18K\log(2d/\delta)}{\xi}$, then it holds with probability at least $1-\delta$ that 
\begin{equation*}
\left(1 - \frac{3\xi + \mu}{3\xi + 3\mu}\right) H^{(t)} \preceq\nabla^{2}F(x^{(t)}) \preceq \left(1+ \frac{L - 2\xi}{2(L+\xi)}\right) H^{(t)}.
\end{equation*}
\end{proof}

\begin{proof}{\bf of Theorem~\ref{thm:NewSamp}}
	Let us denote 
	\begin{equation*}
	H_{\SM} = \frac{1}{|\SM| } \sum_{i\in \SM}\nabla^2f_i(x),\quad\mbox{and}\quad \tH = H_{\SM} + \lambda_{r+1} I,
	\end{equation*}
	where $\lambda_{r+1}$ is the $(r+1)$-th largest eigenvalue of $\nabla^{2}F(x^{(t)})$. 
	By the proof of Theorem~\ref{thm:Reg_subnewton} and Eqn.~\eqref{eq:up_1}, if we choose $|\SM| = \frac{18K\log(d/\delta)}{\lambda_{r+1}}$, then we have
	\begin{equation}
	\label{eq:up_2}
	H_{\SM} \succeq \frac{2}{3} \nabla^{2}F(x^{(t)}) - \frac{\lambda_{r+1}}{3} I.
	\end{equation}
	Moreover, by Eqn.~\eqref{eq:down_1} and choosing $|\SM| = \frac{12K\log(d/\delta)}{\lambda_{r+1}}$, we can obtain that
	\begin{equation}
	\label{eq:down_2}
	H_{\SM} \preceq \frac{3}{2} \nabla^{2}F(x^{(t)}) + \frac{\lambda_{r+1}}{2} I.
	\end{equation}
	By Corollary 7.7.4 (c) of \citet{horn2012matrix}, Eqn.~\eqref{eq:up_2} and \eqref{eq:down_2} imply that
	\begin{equation}
	\label{eq:lam_val}
	\frac{1}{3} \lambda_{r+1} \leq\lambda_{r+1}(H_{\SM}) \leq 2\lambda_{r+1}.
	\end{equation}
Let us express the SVD  of $H_{\SM}^{(t)}$ as follows
\[
H_{\SM}^{(t)} = U\hat{\Lambda}U^T = U_r \hat{\Lambda}_r U^T_r + U_{{\setminus}r} \hat{\Lambda}_{{\setminus}r}U_{{\setminus}r}^T.
\]
Then $H^{(t)}$ can be represented as
\begin{equation*}
H^{(t)} = H_{\SM} +
U
\left[
\begin{array}{cc}
0 & 0\\
0& \lambda_{r+1}(H_{\SM}) I - \hat{\Lambda}_{{\setminus}r}
\end{array}
\right]
U^T.
\end{equation*}
By Eqn.~\eqref{eq:up_2} and $\frac{1}{3} \lambda_{r+1} \leq\lambda_{r+1}(H_{\SM})$ (Eqn.~\eqref{eq:lam_val}), we have
\begin{align*}
H^{(t)} \succeq  \frac{2}{3} \nabla^{2}F(x^{(t)}) - \frac{\lambda_{r+1}}{3} I + U
\left[
\begin{array}{cc}
0 & 0\\
0& \lambda_{r+1}(H_{\SM}) \cdot I - \hat{\Lambda}_{{\setminus}r}
\end{array}
\right]
U^T 
\succeq 
\frac{2}{3} \nabla^{2}F(x^{(t)})
\end{align*}
which implies that
\begin{equation*}
\nabla^{2}F(x^{(t)}) \preceq \left(1+\frac{1}{2}\right)H^{(t)}.
\end{equation*}

By Eqn.~\eqref{eq:down_2} and~\eqref{eq:lam_val}, we have
\begin{align*}
	H^{(t)}
	\preceq& \frac{3}{2} \nabla^{2}F(x^{(t)}) + \frac{\lambda_{r+1}}{2} I
	+
	U
	\left[
	\begin{array}{cc}
	0 & 0\\
	0& \lambda_{r+1}(H_{\SM}) I - \hat{\Lambda}_{{\setminus}r}
	\end{array}
	\right]
	U^T\\
	\preceq&\frac{3}{2} \nabla^{2}F(x^{(t)}) + \frac{5}{2} \lambda_{r+1} I\\
	\preceq& \left(\frac{3}{2} + \frac{5\lambda_{r+1}}{2\mu}\right) \nabla^{2}F(x^{(t)})
\end{align*}
which implies that
\begin{equation*}
\left(1 - \frac{5\lambda_{r+1} + \mu}{5\lambda_{r+1} + 3\mu}\right) H^{(t)} \preceq \nabla^{2}F(x^{(t)}).
\end{equation*}

Therefore, if choosing $|\SM| = \frac{18K\log(2d/\delta)}{\lambda_{r+1}}$, we can obtain that
\begin{equation*}
\left(1 - \frac{5\lambda_{r+1} + \mu}{5\lambda_{r+1} + 3\mu}\right) H^{(t)} \preceq \nabla^{2}F(x^{(t)}) \preceq \left(1+\frac{1}{2}\right)H^{(t)}.
\end{equation*}
The convergence properties can be derived directly by Theorem~\ref{thm:univ_frm}.
\end{proof}

\end{document}